\documentclass[a4paper,11pt]{article}
\usepackage{cancel}
\usepackage{amsmath}
\usepackage{amssymb}
\usepackage{soul}
\usepackage{arydshln}
\usepackage{adjustbox, lipsum}

\usepackage[utf8]{inputenc}
\usepackage{amsfonts}
\usepackage{pstricks}
\usepackage{rotating}
\usepackage{natbib}
\usepackage{booktabs}

\usepackage{authblk}
\usepackage{multirow}

\renewcommand{\arraystretch}{1.2}

\newtheorem{prop}{\bf Proposition}[section]

\newtheorem{rem}{\bf Remark}[section]

\newcommand{\dem}{\par \noindent{\bf Proof:} }
\newcommand{\fin}{\hfill $\square$  \par \bigskip}
\definecolor{gr}{RGB}{0,153,0}

\author[(a)]{Alfredo Marín}
\author[(b)]{Luisa I. Martínez-Merino\footnote{Corresponding author: luisa.martinez@uca.es}}
\author[(b)]{Antonio M. Rodríguez-Chía }
\author[(c,d)]{Francisco Saldanha-da-Gama}

\affil[(a)]{\small{Departamento de Estadística e Investigación Operativa, Facultad de Matemáticas, Universidad de Murcia, Murcia, Spain, amarin@um.es}}
\affil[(b)]{\small{Departamento de Estadística e Investigación Operativa, Universidad de Cádiz, Cádiz, Spain, luisa.martinez@uca.es, antonio.rodriguezchia@uca.es}}
\affil[(c)]{\small{Departamento de Estatística e Investigação Operacional, Faculdade de Ciências da Universidade de Lisboa, Portugal, faconceicao@fc.ul.pt}}
\affil[(d)]{\small{Centro de Matemática, Aplicações Fundamentais e Investigação Operacional, Faculdade de Ciências da Universidade de Lisboa, Portugal, faconceicao@fc.ul.pt}}

\date{}

\title{\LARGE Multi-period Stochastic Covering Location Problems: modeling framework and solution approach}

\allowdisplaybreaks

\begin{document}
\maketitle

\begin{abstract}
	 This paper introduces 
	 a very general discrete covering location model that accounts for uncertainty and time-dependent aspects. 
	 A MILP formulation is proposed for the problem.
	 Afterwards, it is observed that most of the models existing in the literature related with covering location can be considered as particular cases of  this formulation. 
	 In order to tackle large instances 
	 of this problem a Lagrangian relaxation based heuristic is developed. 
	 A computational study is addressed to check the potentials and limits of the formulation and some variants proposed for the problem, as well as to evaluate the heuristic. Finally,
	 different measures to report the relevance of considering a multi-period stochastic setting are studied.

\noindent		
	\textbf{Keywords:} Location, covering, discrete optimization, multi-period facility location, Lagrangian relaxation

\end{abstract}

\section{Introduction}
\label{sec:introduction}
A covering location problem consists of deciding where to locate a set of facilities/equipments in such a way that a set of demand points can be reached or served within some previously defined time or distance; in that case, a demand point is referred to as covered.
When the set of potential locations for the facilities is discrete, we are facing a discrete covering location problem.
This is the case we focus in this work. 
{Some reviews of covering location problems can be found in \cite{Schilling93}, \cite{Berman2010}, \cite{Snyder2011}, \cite{Farahani2012}, and \cite{Garcia2015}.}

Like for many other problems, the use of optimization techniques in the context of covering facility location calls for an \textsl{a priori} development of mathematical models.
In this respect,
two main modeling frameworks can be found in the literature depending on the 
objective of the problem: the set covering location problem (SCP) and the maximal covering location problem (MCLP).
The former was initially introduced by \cite{Toregas71} and aims at minimizing the total cost associated with the installed facilities that guarantees the coverage of all demand points. 
The latter was initially proposed by \cite{Church74} and consists of locating a certain number of facilities maximizing the resulting coverage. 

One interesting and relevant feature of covering location problems is that they can be applied to different areas.
Such applications
often regard the location of emergency or surveillance facilities such as ambulances, fire fighting bases, or surveillance points. 
In such a case, by covering the ``demand points'' with a set of facilities we are ensuring that a certain response level is guaranteed upon the occurrence of an event.
The SCP has been applied to deployment of emergency services \citep{Toregas71} and to the analysis of 
markets \citep{Storbeck88}.
On the other hand, the MCLP has been 
applied, for instance, to cluster and discriminant analysis \citep{Chung86}.
Overall, the flexibility provided by covering location problems explains why they have been deeply studied. 
For further details, the reader can refer to the reviews by \cite{Schilling93} and \cite{Caprara2000} and to the references therein.

In this paper we aim at bringing further the comprehensiveness of the optimization models in use when it comes to discrete covering location problems. 
With this purpose, we extend the general modeling framework proposed by \cite{Garcia2015} in order to capture two features (alone or combined) of practical relevance: time-dependency of aspects such as demand and costs in addition to the uncertainty that often underlies the parameters of the problem.
It is worth noticing that the model proposed by \cite{Garcia2015} is already very general since it unifies, among others, the set covering location problem and the maximal covering location problem.
Nevertheless, the extension we propose in this 
paper has that unified version as a particular case.
We note that it is not our purpose in this work
to study a specific problem that can be directly applied to some practical situation (although we do not deny that either).
Instead, we are investigating a broad family of problems under a common ``umbrella'', i.e., we focus on a general problem that includes as particular cases many problems of practical relevance. 
Accordingly, all the developments we propose 
are directly applicable to all those particular cases. 

We consider 
a finite planning horizon that is partitioned into several time-periods (not necessarily of the same length). 
In each period, multiple facilities/equipments can be opened or closed in each 
location at the expenses of the corresponding 
costs.
We assume that in each period it is desirable to cover 
each customer by at least a certain number of facilities/equipments.
If coverage goes beyond the minimum threshold, some benefit is achieved; otherwise we have a penalty for uncovered  demand. 
Both the benefits and penalties are assumed to be time-dependent.
We also consider the case in which the demand points covered by each potential location may change over time. 
This is relevant since the covering capability of some facilities may depend not only on the distance but also on other 
features such as the traveling time (it is often not the same in the winter with snow on the roads as in the summer with dry roads), visibility (e.g., in a forest across the year), aging of equipment rendering
less coverage over time, \textit{et cetera}. 

The second relevant feature that we capture in this work is uncertainty.
We 
assume that uncertainty is associated with the demand for coverage and also with the capability of a facility to cover the demand points.
Moreover, we assume that uncertainty can be 
captured by a random vector with a finite support that has a
joint probability distribution {function} known in advance.
If the support of the random vector is not finite, 
we may have to resort to some sampling (and thus approximate) method that nonetheless will require the resolution of a problem such as the one we are considering.
Each realization of the underlying random vector will be referred to as a \textsl{scenario}.
In other words, a scenario fully specifies all the {uncertain parameters.} 

With the purpose of gathering time and uncertainty as additional dimensions in covering facility location we propose a two-stage stochastic programming model.
The \textit{ex ante} or \textit{first-stage} decisions are associated with the definition of a location plan to be implemented throughout the planning horizon. 
The \textit{ex post} or \textit{second-stage} decisions are dependent on how uncertainty is ``revealed'' and include the definition of a coverage level for each demand point that, in turn, allows the evaluation of the corresponding costs.
The goal is to minimize the total expected cost, which includes the cost associated with the facilities (opening/closing/operation) and the total expected penalty for coverage shortage minus the total benefit for coverage surplus.

As our computational experiments show, hardly can the model we propose be tackled by a general-purpose solver when it comes to solving medium- or large sized instances. 
For that reason, in addition to a contribution
in terms of the mathematical model proposed, we develop a Lagrangian relaxation based heuristic.
We handle the Lagrangian dual using subgradient optimization.
In each iteration, a feasible solution for the problem is derived. 

When it comes to capturing features such as the existence  of time-dependent parameters or uncertainty, one relevant question is related with the advantages of considering a more complex optimization model instead of a simplified one (e.g., static or deterministic).
An answer can be given by two distinct measures: the Expected Value of Perfect Information (EVPI), and the Value of Multi-period Solution (VMS).
The first (see for instance, {\citealt{Correia2015}} for a discussion in the context of facility location problems) represents the maximum amount that the decision maker is willing to pay to get perfect information about the future. 
The larger the observed value 
the more relevant may be to deal with uncertainty. 
The latter measure was first introduced by {{ \cite{Alumur&Nickel&Saldanha-da-Gama&Verter:2012}}} and formalized afterwards by { \cite{Nickel2015}}. 
It looks 
into an approximate static solution and evaluates it in the multi-period model. 
By doing so, we get an indication in terms of the relevance of explicitly considering a multi-period setting.

As it was mentioned above, the model introduced in this work takes as a starting point the general model proposed by \cite{Garcia2015}. 
We refer the reader to that book chapter for an overview in terms of the existing literature on deterministic and static covering facility location problems.

Multi-period problems in the context of covering facility location have been recently revisited by \cite{Seyed2016}, who included a specific section on dynamic covering location problems.
The authors describe the problems and optimization models existing in the literature such as those introduced 
by \cite{Gun82} and that will also be generalized in the current work.
Additionally, the above-mentioned review paper quotes the works by \cite{Gendreau2001} and \cite{Rajagopalan2008}.
The former proposes a dynamic covering location model for real time ambulance relocation. 
The objective is to maximize the proportion of demand covered by at least two vehicles within a certain radius minus a relocation cost. 
The latter focuses on ambulance allocation in certain time clusters to meet some coverage requirements.

Another type of problem in the context of covering location problems is the one proposed by \cite{Colombo2016}, which maximizes the demand coverage given a set of demand points that require different types of service.
Before that, \cite{Zarandi2013} had introduced a multi-period MCLP.
In that problem the number of centers to be opened throughout the entire planning horizon is restricted. 
 
It must be noted that none of the works cited above address uncertainty in the context of multi-period covering location problems. 
To the best of the authors' knowledge, the only work doing so is due to {\cite{Vatsa2016}}, who investigate a multi-period maximal covering model with server uncertainty. 
The authors adopt a robust optimization approach. 
In our case, as explained before, we are adopting a stochastic programming modeling framework, which allows us to explicitly embed a probability distribution in our models as a way for capturing uncertainty.

{In summary, the contribution of this paper is to provide a very general model of a covering location problem capturing uncertainty and time-dependent aspects. This model includes as particular cases many of the models studied in the literature. In addition, a specific Lagrangian relaxation based heuristic has been developed for solving medium or large sized instances of this problem.}

The remainder of this paper is organized as follows. 
In Section \ref{general_model}, we provide all the details related with our general problem and we introduce an optimization model for it. 
Section \ref{sec:Lagrangian-heuristic} is devoted to the development of a Lagrangian relaxation based heuristic that provides both lower and upper bounds on the optimal solutions.
The numerical experiments conducted with the models and methodologies proposed in this paper are reported and discussed in Section \ref{sec:numerical_experiments}. 
Finally, the paper ends with an overview of the work done.

\section{Modeling framework}
\label{general_model}
We consider a general multi-period stochastic covering facility location problem.
The problem consists of making a decision for a planning horizon of finite length that is divided into several time periods. 
A set of potential locations for the facilities is given as well as a set of demand points that must  be desirably covered throughout the planning horizon.
In order to consider a more general setting, we assume that in each location we can install more than one facility/equipment.
However we consider a limit for the number of facilities in each location and we assume that in each time period there is a limit for the total number of facilities installed across all locations.

It is often the case that we need to make decisions for systems that are already operating.
We consider that this may also be the case and thus we assume that in some locations we may have some facilities currently operating.

Opening or closing a facility during the planning horizon will suppose the payment of some costs. Besides, we assume that the opening of a facility always occurs at the beginning of a time period and the closing always occurs at the end of the period. In addition to that, we consider that operating a facility also incurs a cost. 
In what follows, we will claim that a facility is opened at the beginning of {some period} $t$ if this facility is not operating in time period $t-1$, but it operates {in 
$t$}. 
Equivalently, we will claim that a facility is closed at the end of period $t$ if it is operating in $t$ but it is not in $t+1$.

We assume that uncertainty is associated both with the {covering capability of each pair $(i,j)$, with $i$ denoting a potential facility location and $j$ a demand point,}
and with the minimum threshold for the coverage of each demand point in each time period.
We consider unit-demand customers, that is, our goal is simply to define a plan for covering the customers in the best possible way.
We suppose that a benefit (that is also uncertain) is associated with coverage surplus (i.e. the number of facilities/equipments covering a demand point above a minimum threshold) and a penalty (again subject to uncertainty) is associated with the coverage shortage (i.e. the number of facilities/equipments in shortage for a demand point with respect to the minimum threshold).

Finally, we consider that uncertainty can be captured by a finite set of scenarios each having a probability assumed known in advance.
{Moreover,}
we consider that one scenario fully specifies all the stochastic quantities.

The decisions to make are twofold: the first-stage decisions that include the opening/closing plan for all locations throughout the planning horizon and the second-stage decisions that help accounting for the total benefit and penalty achieved by the covering plan.

Next, we introduce the relevant notation to be used hereafter.

\medskip
\noindent
Sets:

\begin{tabular}{p{30mm}p{115mm}}
	${\cal T}= \{1,\ldots,T\}$, & set of periods in the planning horizon.
\end{tabular}

\begin{tabular}{p{30mm}p{115mm}}
	$I= \{1,\ldots,m\}$, & set of potential location for the facilities.
\end{tabular}

\begin{tabular}{p{30mm}p{115mm}}
	$J=\{1,\ldots,n\}$, & set of demand points.
\end{tabular}

\begin{tabular}{p{30mm}p{115mm}}
	${\cal S}= \{1,\ldots,S\}$, & set of scenarios.
\end{tabular}

\medskip
\noindent
Deterministic costs:

\begin{tabular}{p{5mm}p{135mm}}
	$o_{it}$, & cost of opening a facility at $i \in I$ at the beginning of time period $t \in {\cal T}$.
\end{tabular}

\begin{tabular}{p{5mm}p{135mm}}
	$c_{it}$, & cost of closing a facility at $i \in I$ at the end of time period $t \in {\cal T} \setminus\{ T\} $.
\end{tabular}

\begin{tabular}{p{5mm}p{135mm}}
	$f_{it}$, & cost  of operating a facility at $i 	 \in I$ during period $t \in {\cal T}$.
\end{tabular}

\medskip
\noindent
Other deterministic parameters:

\begin{tabular}{p{5mm}p{135mm}}
	$e_i$, & maximum number of facilities/equipments  that can be operating in location $i \in I$.
\end{tabular}

\begin{tabular}{p{5mm}p{135mm}}
	$p_t$, & maximum number of facilities/equipments that can be  operating in time period $t \in {\cal T}$.
\end{tabular}

\begin{tabular}{p{5mm}p{135mm}}
	 $ \bar{y}_{i0}$, & number of facilities that are open at location $i \in I$ before the beginning of the planning horizon.
\end{tabular}

\medskip
\noindent
Stochastic (scenario-dependent) parameters

\begin{tabular}{p{5mm}p{135mm}}
	$a_{ijt}^s$, & (binary) parameter indicating whether location $i\in I$ can cover demand point  $j\in J$ in time period  $t\in {\cal T}$ under scenario $s \in {\cal S}$.
\end{tabular}

\begin{tabular}{p{5mm}p{135mm}}
	$b_{jt}^s$, & minimum  number of facilities/equipments requested to cover demand point $j \in J$ in period $t \in {\cal T}$ under scenario $s \in {\cal S}$.
\end{tabular}

\begin{tabular}{p{5mm}p{135mm}}
	$d_{jkt}^s$, & marginal benefit for covering demand point $j \in J$ with a surplus of at least $k$ facilities in period $t \in {\cal T}$ under scenario $s \in {\cal S}$.
\end{tabular}

\begin{tabular}{p{5mm}p{135mm}}
	$h_{jkt}^s$, & marginal penalty for a shortage of at least $k$ facilities in the coverage of demand point $j \in J$ in period $t \in {\cal T}$ under scenario $s \in {\cal S}$.
\end{tabular}

\begin{tabular}{p{5mm}p{135mm}}
	$\pi_s$,  & probability that scenario $s \in {\cal S}$ occurs. We assume that  $\pi_s>0$ for $s\in {\cal S}$ and $\sum_{s \in {\cal S}} \pi_s = 1$.
\end{tabular}

\medskip
\noindent
Auxiliary sets:

\medskip
\begin{tabular}{p{38mm}p{100mm}}
	$K_{jt}^s = \{1,\ldots,p_t - b_{jt}^s\}$,  & with $p_t - b_{jt}^s$ representing the maximum surplus in terms of coverage that we may observe in demand point $j \in J$ in period $t \in {\cal T}$ under scenario $s \in {\cal S}$.
\end{tabular}

\begin{tabular}{p{38mm}p{100mm}}
	$K_{jt}^{\prime s} = \{1,\ldots,b_{jt}^s\}$, & with $b_{jt}^s$ standing for the maximum shortage in terms of coverage that we may observe in demand point $j \in J$ in period $t \in {\cal T}$ under scenario $s \in {\cal S}$. 
\end{tabular}

\bigskip
In the model proposed next, we are accounting for the total (expected) cost. 
Accordingly, the benefit for coverage surplus will be treated as a negative cost. 
For this reason, in addition to the above notation we consider

\medskip
\begin{tabular}{p{25mm}p{118mm}}
	$g_{jkt}^s = -d_{jkt}^s$, & the negative marginal cost representing the (marginal) benefit accounted for by $d_{jkt}^s$.
\end{tabular}

\medskip
In our work, the costs $g$ and $h$ are being looked at as marginal variations in the cost for having one extra facility in surplus or shortage. 
For instance, suppose that at some demand point $j$ in some period $t$ and scenario $s$ we have a surplus of $2$ facilities/equipments. 
In this case, we assume the corresponding ``cost'' to be $g_{j1t}^s + g_{j2t}^s$. 
This means that $g_{j2t}^s$ can be looked at as the increase in the ``cost'' for going from one to two surplus facilities.
If the surplus becomes $3$ then the ``cost'' would be $g_{j1t}^s + g_{j2t}^s + g_{j3t}^s$. 
Hence, an increase of one surplus facility when we have already two surplus facilities renders an increase in the benefit of $g_{j3t}^s$.

We assume that for each  demand point $j\in J$, each period $t\in {\cal T}$ and scenario $s \in {\cal S}$ $$g_{j1t}^s \leq g_{j2t}^s \leq \ldots \leq g_{j|K_{jt}^s|t}^s  \leq 0$$
 and 
 $$0  \leq h_{j1t}^s\leq h_{j2t}^s\leq\ldots\leq h_{j|K_{jt}^{\prime s}|t}^s.$$ 

In other words, we are assuming non-increasing marginal benefit (for surplus) and a non-decreasing marginal penalty (for shortage).
We believe 
that this is the case that holds in practice: the marginal benefit of having a larger surplus decreases when the number of surplus facilities increases; the marginal penalty for having shortage of facilities increases with the number of facilities in shortage.

An integer programming model can be presented for our problem considering the following five sets of decision variables:
\begin{eqnarray*}
	z_{it} & = & \textrm{Number of facilities opened in site $i \in I$ at the beginning of period $t \in {\cal T}$}. \\
	z^\prime_{it} & =&\textrm{Number of facilities closed in location $i \in I$ at the end of period $t \in {\cal T} \setminus \{ T\}$}. \\
	y_{it}&=&\textrm{Number of facilities operating in location $i \in I$ in period $t \in {\cal T}$}. \\	
	w_{jkt}^s & = &
	\begin{cases}
		1, &\textrm{if demand point $j$ is covered by at least $ b_{jt}^s+k$ facilities in period $t$}\\
			& \textrm{ and scenario $s$,}\\
		0, &\textrm{otherwise},
	\end{cases}\\
	&&\textrm{with }j \in J,\;k\in K_{jt}^s,\;t\in {\cal T},\;s \in {\cal S}.\\
	 v_{jkt}^s & = &
	\begin{cases}
		1, &\textrm{if demand point $j$ is covered by at most $b_{jt}^s-k$ facilities in period $t$}\\
		& \textrm{ and scenario $s$,}\\
		0, &\textrm{otherwise},
	\end{cases}\\
	&&\textrm{with }j\in J,\;k\in K_{jt}^{\prime s},\;t\in {\cal T},\; s \in {\cal S}.
\end{eqnarray*}

Finally, we can write the following integer programming model that we denote by \linebreak (GMSCLP), i.e, General Multi-period Stochastic Covering Location Problem:
\begin{alignat}{3}
& \mbox{min} & \quad & \sum_{i\in I}\sum_{t\in {\cal T}}o_{it}z_{it}+\sum_{i\in I}\sum_{t\in {\cal T}\setminus\{ T\}} c_{it}z^\prime_{it}+\sum_{i\in I}\sum_{t\in {\cal T}}f_{it} y_{it} & & \hspace{25ex} \nonumber \\
& & & + \sum_{s\in {\cal S}}\pi_s \left[\sum_{t\in {\cal T}}\sum_{j\in J}\sum_{k\in K_{jt}^s}g_{jkt}^s w_{jkt}^s+\sum_{t\in {\cal T}}\sum_{j\in J}\sum_{k\in K_{jt}^{\prime s}}h_{jkt}^sv_{jkt}^s \right] \label{obj}
\end{alignat}
\vspace{-5mm}
\begin{alignat}{3}
& \mbox{s.t.} & \quad & \sum_{i\in I}y_{it}\leq p_{t}, & & t\in {\cal T}, \label{sumy} \\
& & & y_{i1}=z_{i1}+ \bar{y}_{i0},  & & i\in I, \label{yz}\\
& & & y_{it}=y_{i,t-1}+z_{it}-z^\prime_{i,t-1}, & & i\in I,t\in {\cal T}\setminus\{1\}, \label{yz2}\\
& & & \sum_{i\in I} a_{ijt}^sy_{it}=b_{jt}^s+\sum_{k\in K_{jt}^s}w_{jkt}^s-\sum_{k\in K_{jt}^{'s}}v_{jkt}^s, & \quad & j\in J,t\in {\cal T},s\in {\cal S},\label{cover}\\
& & & w_{j1t}^s+v_{j1t}^s\leq 1, & & j\in J,t\in {\cal T},s\in {\cal S},\label{ww1}\\
& & & w_{jkt}^s \leq w_{j1t}^s, & & j\in J,t\in {\cal T},s\in {\cal S},k=2,\ldots,|K_{jt}^s|,\label{ww2}\\
& & & v_{jkt}^s \leq v_{j1t}^s, & & j\in J,t\in {\cal T},s\in {\cal S},k=2,\ldots,|K_{jt}^{'s}|,\label{ww3}\\
& & & y_{it}\in\{0,\ldots,e_i\}, & & i\in I,t\in {\cal T}, \label{ybound}\\
& & & z_{it}\in\{0,\ldots,e_i\}, & & i\in I,t\in {\cal T}, \label{zbound}\\
& & & z^\prime_{it}\in\{0,\ldots,e_i\}, & & i\in I,t\in {\cal T}\setminus\{T\}, \label{zpbound} \\
& & & w_{jkt}^s\in\{0,1\}, & & j\in J,k\in K_{jt}^s,t\in {\cal T},s\in {\cal S}, \label{wbound} \\
& & & v_{jkt}^s\in\{0,1\}, & & j\in J,k\in K_{jt}^{\prime s},t\in {\cal T},s\in {\cal S}. \label{vbound}
\end{alignat}

In this model, the objective function \eqref{obj} quantifies the total cost for opening, closing and operating the facilities (first line) plus the total expected benefit plus penalty associated with the actual coverage of the demand points.
Constraints \eqref{sumy} restrict the number of operating centers to $p_t$, for each $t\in {\cal T}$ and constraints \eqref{yz}--\eqref{yz2} establish the 
relationship between $y$~-, $z$- and $z^\prime$-variables.
From a stochastic programming 
point of view, the above constraints together with the domain constraints \eqref{ybound}--\eqref{zpbound} are actually the so-called first-stage constraints in the sense that they only relate the first-stage decision variables.
The left hand side of constraints \eqref{cover} accounts the number of open facilities 
that cover a given customer $j$ in a given scenario $s$ and period $t$; then the right hand side
of the equality aims to compare this amount with the minimum request $b_{jt}^s$, in order to 
ensure that both the $v$- and $w$-variables take their correct values. 
Additional constraints \eqref{ww1}--\eqref{ww3} are 
required to make either all the $v$-variables or all the $w$-variables (or all of them) 
take value 0. The fact that the costs $g_{jkt}^s$ (respectively $h_{jkt}^s$) are 
increasingly sorted guarantees the correctness of the solution.
Finally, constraints \eqref{ybound}--\eqref{vbound} define the domain of the decision variables.

{Using the standard terminology in Stochastic Programming (see, e.g.,  \citealt{Birge&Louveaux:2011}) the above model is the so-called extensive form of the deterministic equivalent since {we} are explicitly considering a (deterministic) model that gathers the first-stage decision and the second stage decision for every possible scenario (the scenarios are extensively considered).
In practice, only one second stage decision (out of all found by the model) is implemented---the one corresponding to the scenario eventually observed.}

A close look into the decision variables shows that we can write the $y$-variables as function of the $z$- and $z^\prime$-variables.
In fact, we have:
\begin{equation*}
y_{it}= \bar{y}_{i0}+\sum_{\tau=1}^{t}z_{i\tau}-\sum_{\tau=1}^{t-1}z^\prime_{i\tau}, \qquad i \in I,\: t \in {\cal T}.
\end{equation*}
Accordingly, we can derive the following alternative model that we denote by (GMSCLP$^\prime$):
\vspace{-0.1cm}
\begin{alignat}{3}
& \mbox{min} & \quad & \sum_{i\in I}\sum_{t\in {\cal T}}o_{it}z_{it}+
\sum_{i\in I}\sum_{t\in {\cal T}\setminus\{T\}}c_{it}z^\prime_{it}+\sum_{i\in I}\sum_{t\in {\cal T}}f_{it}\left[ \bar{y}_{i0}+\sum_{\tau=1}^{t}z_{i\tau}-\sum_{\tau=1}^{t-1}z^\prime_{i\tau}\right]  & & \nonumber \\
& & & +\sum_{s\in {\cal S}}\pi_s \left( \sum_{t\in {\cal T}} \sum_{j\in J} \sum_{k\in K_{jt}^s} g_{jkt}^s w_{jkt}^s+\sum_{t\in {\cal T}}\sum_{j\in J}\sum_{k\in K_{jt}^{'s}}h_{jkt}^s v_{jkt}^s\right) \label{model2-obj}
\end{alignat}
\vspace{-0.5cm}
\begin{alignat}{3}
& \mbox{s.t.} & \quad &  \eqref{ww1}-\eqref{ww3},\eqref{zbound}-\eqref{vbound}, & \qquad & \hspace{4cm}\nonumber \\
& & & \sum_{i\in I}\left(\sum_{\tau=1}^{t}z_{i\tau}-\sum_{\tau=1}^{t-1}z^\prime_{i\tau}\right)\leq p_t-\sum_{i\in I} \bar{y}_{i0}, & & t\in {\cal T},\label{ptrest} \\
& & & \sum_{i\in I}a_{ijt}^s\left( \bar{y}_{i0} + \sum_{\tau=1}^{t}z_{i\tau} - \sum_{\tau=1}^{t-1} z^\prime_{i\tau}\right) =  & & \nonumber \\
& & & \hspace{2cm} = b_{jt}^s + \sum_{k\in K_{jt}^s}w_{jkt}^s - \sum_{k\in K_{jt}^{\prime s}}v_{jkt}^s, & & j\in J,t\in {\cal T},s\in {\cal S},\label{covering}\\
& & &  \bar{y}_{i0}+\sum_{\tau=1}^{t}z_{i\tau}-\sum_{\tau=1}^{t-1}z^\prime_{i\tau}\leq e_i, & & i\in I,\;t\in {\cal T},\label{ytub} \\
& & &  \bar{y}_{i0}+\sum_{\tau=1}^{t}z_{i\tau}-\sum_{\tau=1}^{t-1}z^\prime_{i\tau}\geq 0, & & i\in I,\;t\in {\cal T}.\label{ytlb}
\end{alignat}

Equivalently, the objective function can be expressed as:
\begin{multline}
\sum_{i\in I}\sum_{t\in {\cal T}}\left(o_{it}+\sum_{\tau=t}^{ T}f_{i\tau}\right)z_{it}+
\sum_{i\in I}\sum_{t\in {\cal T}\setminus\{ T\}}\left(c_{it}-\sum_{\tau=t+1}^{ T}f_{i\tau}\right)z^\prime_{it}+
\sum_{i\in I}\sum_{t\in {\cal T}}f_{it}{ \bar{y}_{i0}} \\
+ \sum_{s\in {\cal S}}\pi_s\left(\sum_{t\in {\cal T}}\sum_{j\in J}\sum_{k\in K_{jt}^s}g_{jkt}^s w_{jkt}^s+\sum_{t\in {\cal T}}\sum_{j\in J}\sum_{k\in K_{jt}^{\prime s}}h_{jkt}^sv_{jkt}^s\right).\nonumber
\end{multline}

\begin{rem}\label{remark1}
It is worth noticing that constraints \eqref{ww1}--\eqref{ww3} can be replaced  by other sets of constraints, such as
\begin{eqnarray}
\sum_{k\in K_{jt}^s}w_{jkt}^s&\leq& (1-v_{j1t}^s)(p_t-b_{jt}^s),\ \ \ \,j\in J,t\in {\cal T},s\in {\cal S},\label{opt21}\\
v_{jkt}^s&\leq&v_{j1t}^s,\,\ \ \ j\in J,t\in {\cal T},s\in {\cal S},{k=2,\ldots,|K_{jt}^{\prime s}|}.\label{opt22}
\end{eqnarray}
\noindent or alternatively
\begin{eqnarray}
\sum_{k\in K_{jt}^{'s}}v_{jkt}^s&\leq&(1-w_{j1t}^s)b_{jt}^s,\ \ \ j\in J,t\in {\cal T},s\in {\cal S},\label{opt31}\\
w_{jkt}^s&\leq&w_{j1t}^s,\ \ \ j\in J,t\in {\cal T},s\in {\cal S}, {k=2,\ldots,|K_{jt}^{s}|} .\label{opt32}
\end{eqnarray}	

These two sets of constraints are weaker than \eqref{ww1}--\eqref{ww3}. 
{In fact,} {\eqref{opt21} can be obtained from the aggregated version of constraints \eqref{ww2} by summing up $k$ and, by \eqref{ww1}, replacing the resulting} {right-hand} {side, i.e., $(p_t-b_{jt}^s)w_{j1t}^s$ by $(p_t-b_{jt}^s)(1-v_{j1t}^s)$. Analogously we can argue for \eqref{opt31}. Moreover, the two sets of constraints together with \eqref{ww2} and \eqref{ww3} respectively, analogously to \eqref{ww1}}{--}{ \eqref{ww3}, make either the $v$-variables or the $w$-variables (or all of them) take value $0$.}
However, the number of constraints they involve is smaller.
We have addressed preliminary computations showing that constraints \eqref{opt21}--\eqref{opt22} and \eqref{opt31}--\eqref{opt32} are ineffective regarding CPU times and, as expected, {in terms of the} linear relaxation bounds.
\begin{flushright}
$\tiny \blacksquare$
\end{flushright}
\end{rem}

\subsection{The particular cases}
\label{subsec:particularCases}
Although an extensive literature can be found on covering location problems, not many works focus a multi-period or stochastic nature that some of these problems may have.
In this section we describe how the model (GMSCLP) generalizes most of the models presented in the literature, which shows the generality of the modeling framework we have just proposed. 

Since model (GMSCLP) was built upon the deterministic and static (single-period) case (COV) proposed by \cite{Garcia2015}, it immediately generalizes that model as well as all the particular cases found in that book chapter, such as the SCP, the weighted SCP, the redundant covering location problem, \textit{et cetera}. 

Model (GMSCLP) reduces to the model (COV) from \cite{Garcia2015} if we set $T=\{1\}$, $S=\{1\}$, $\pi_1=1$, and $ \bar{y}_{i0}=0$ for $i\in I$. 
In this case, $z$- and $z^\prime$- variables are not necessary; only $y$-variables are needed to decide whether some facilities are opened in location $i$ or not. 
Specifically, the parameters that adapt our formulation to (COV) are: 

$o_{i1}=0$, $c_{i1}=0$, $f_{i1}=f_i$, for $i\in I$; $p_1=p$;

$b_{j1}^1=b_j$, $K_{j1}^1=K_j$; 

$g_{jk1}^1=g_{jk}$ for $j\in J,k\in K_{j1}^1=K_j$ for $j\in J$; and

$K_{j1}^{'1}=\emptyset$ for $j\in J$.

\smallskip
Additionally, we note that model (COV) does not capture the case in which a demand point $j$ is covered by less than $b_j$ facilities.
Consequently, we have that $K^{\prime}_j=\emptyset$ for $j\in J$ and $h$ vector is not considered in this model.

Concerning existing multi-period covering facility location problems, we have already quoted the paper by \cite{Gun82}, where different problems and models are investigated.
Our general formulation also generalizes the models proposed in that work as we 
show next.

Since all the problems studied in the above-mentioned paper
are deterministic, they can be looked at as single-scenario problems, i.e, we have 
$S=\{1\}$ with $\pi_1=1$. 
Thus, for the sake of simplicity in terms of notation, in the following particularizations, index $s$ is 
omitted.

To start with, it is important to note that if we take $K_{jt}=K^{^\prime}_{jt}=\emptyset$, then we can replace the  equalities in constraints \eqref{cover} by inequalities yielding
\begin{equation}
\sum_{i\in I} a_{ijt}y_{it}\geq b_{jt},\;j\in J,t\in {\cal T}.\label{ineqcov}
\end{equation}

In some of the models proposed by \cite{Gun82}, facilities cannot be closed during the planning horizon. 
This can be handled in our model by assigning
a very large value ($M$) 
to the closing costs ($c_{it}$). 
In this way, we ensure that facilities can only be opened 
during the planning horizon. 
Similarly, when centers cannot be opened, the same value $M$ is assigned to opening costs ($o_{it}$).

{
Next, we detail 
how our model can be particularized to get the most general multi-period models proposed by \cite{Gun82}.

\begin{itemize}
   \item GDSCLP---General Dynamic Set Covering Location Problem.
	
	For this problem, in some (empty) potential locations ($I^o\subset I$) we can open facilities during the planning horizon; in the other locations $ (I^{c}\subset I )$ we have a facility operating that can be closed during the planning horizon. 
	These two sets do not overlap.
	For each $t\in {\cal T}$, a subset of customers $J_t\subset J$ must be covered.
	
	To adapt model (GMSCLP) we can remove $w$- and $v$-variables.
	Constraints \eqref{cover} are replaced by constraints \eqref{ineqcov}. 
	For the remaining parameters we set: $K_{jt}=K'_{jt}=\emptyset$, $f_{it}=1$, $p_t=m$, $e_i=1$, for $i\in I$, $j\in J$, $t\in {\cal T}$. 
	Besides, $o_{it}=M$ for $i \in I^c$ and $o_{it}$=0 otherwise for $t\in \cal T$. Similarly, $c_{it}=M$ for $i\in I^o$ and $c_{it}=0$ otherwise.
	Finally, $b_{jt}=1$ for $j\in J_t$, $t\in {\cal T}$ and $b_{jt}=0$ for $j\notin J_t$, $t\in {\cal T}$.
	
	\item DMCLP2---Dynamic Maximum Covering Location Problem 2.
	
	This problem is a multi-period maximum covering location problem which minimizes the weighted sum of demand points that are not covered in each time period plus the opening and closing costs in each period. 
	Each weight $(\gamma_{jt})$ is the cost of demand point $j$ when it is not covered in period $t$. 
	In this case, the facilities can be opened or closed during the planning horizon. 
	
	To model this problem using our general formulation we set $K_{jt}=\{1\}$ and  $g_{j1t}=\gamma_{jt}$ for $j\in J$, $t\in {\cal T}$.
	Constraints \eqref{cover} are replaced by
	\begin{equation}
	\sum_{i\in I} a_{ijt}y_{it}\geq w_{j1t},\ \ \ j\in J,t\in {\cal T}. \label{ineqcov2}
	\end{equation}
	For the remaining parameters, $o_{it}=o_t$, $c_{it}=c_t$, $f_{it}=0$, $K^\prime_{jt}=\emptyset$, $p_t=p_t$, $b_{jt}=0$ and $e_i=1$ for $i\in I$, $j\in J$, $t\in {\cal T}$. 
	To obtain the same objective value as in the original model, constant $-\sum_{j\in J}\sum_{t\in {\cal T}}g_{j1t}$ must be added to the objective value.	
\end{itemize}

In tables \ref{const} and \ref{param} of the Appendix we summarize the constraints and parameters, respectively, used to adapt (GMSCLP) to the models described above. 
We also include in {the}
Appendix the analysis of the remaining models proposed in \cite{Gun82}.
}

The general problem we are investigating in this work is NP-hard since it contains as special cases some well-known NP-hard problems (e.g. the MCLP---see { \citealt{Megiddo&Zemel&Hakimi:1983})}.
As our computational experiments reported in Section~\ref{sec:numerical_experiments} show, only for small instances of the problem can we expect a general-purpose solver to be effective for tackling them.
This motivates the developments we propose next.

\section{A Lagrangian relaxation based procedure}
\label{sec:Lagrangian-heuristic}
In this section we develop a Lagrangian relaxation based algorithm for finding (hopefully) high quality feasible solutions to the general problem introduced in the previous section.

In the methodological developments proposed in this section, we consider model (GMSCLP$^\prime$).
In particular, we note that $y$-variables will not be used. 

We start by relaxing the covering constraints \eqref{covering} in a Lagrangian fashion by considering multipliers $\alpha_{jt}^s\in \mathbb{R}$. 
By doing so, we obtain the following relaxed problem, defined for each vector ${\bf \alpha}$ that we denote by ($\mbox{LR}_\alpha$), 
\begin{eqnarray}
&\min&
\sum_{i\in I}\sum_{t\in {\cal T}}\left(o_{it}+\sum_{\tau=t}^{ T}f_{i\tau}\right)z_{it}+
\sum_{i\in I}\sum_{t\in {\cal T}\setminus\{ T\}}\left(c_{it}-\sum_{\tau=t+1}^{ T}f_{i\tau}\right)z'_{it}+\nonumber\\
&&\sum_{i\in I}\sum_{t\in {\cal T}}f_{it} \bar{y}_{i0}+\sum_{s\in {\cal S}}\pi_s\left(\sum_{t\in {\cal T}}\sum_{j\in J}\sum_{k\in K_{jt}^s}g_{jkt}^s w_{jkt}^s+\sum_{t\in {\cal T}}\sum_{j\in J}\sum_{k\in K_{jt}^{'s}}h_{jkt}^sv_{jkt}^s\right)+\nonumber\\
&&\sum_{j\in J}\sum_{t\in {\cal T}}\sum_{s\in {\cal S}}\alpha_{jt}^s\left(\sum_{i\in I}a_{ijt}^s\left( \bar{y}_{i0}+\sum_{\tau=1}^{t}z_{i\tau}-\sum_{\tau=1}^{t-1}z'_{i\tau}\right)\right)+\nonumber\\
&&\sum_{j\in J}\sum_{t\in {\cal T}}\sum_{s\in {\cal S}}\alpha_{jt}^s\left(-b_{jt}^s-\sum_{k\in K_{jt}^s}w_{jkt}^s+\sum_{k\in K_{jt}^{'s}}v_{jkt}^s\right)\hspace{1.5cm}\nonumber\\
& &\nonumber \\
& \mbox{s.t.} &  \eqref{ww1}-\eqref{ww3}, \eqref{zbound}-\eqref{vbound}, \eqref{ptrest}, \eqref{ytub},\eqref{ytlb}.\nonumber
\end{eqnarray}

In turn, this model can be split in two subproblems, denoted by ($\mbox{LR1}_\alpha$) and  ($\mbox{LR2}_\alpha$) as follows:
\begin{eqnarray*}
 \textrm{($\mbox{LR1}_\alpha$)}	&\min&\sum_{i\in I}\sum_{t\in {\cal T}}\left[o_{it}+\sum_{\tau=t}^{ T}(f_{i\tau}+\sum_{j\in J}\sum_{s\in {\cal S}}\alpha_{j\tau}^s a_{ij\tau}^s)\right]z_{it}+\\
	&&\sum_{i\in I}\sum_{ t\in {\cal T}\setminus \{T\}}\left[c_{it}-\sum_{\tau=t+1}^{T}(f_{i\tau}+\sum_{j\in J}\sum_{s\in {\cal S}}\alpha_{j\tau}^s a_{ij\tau}^s)\right]z'_{it}+\sum_{i\in I}\sum_{t\in {\cal T}}(f_{it}+\sum_{j\in J}\sum_{s\in {\cal S}}\alpha_{jt}^s a_{ijt}^s) \bar{y}_{i0} \\
	& & \\
	& \mbox{s.t.} & \eqref{zbound},\eqref{zpbound},
		\eqref{ptrest},\eqref{ytub},\eqref{ytlb}.
\end{eqnarray*}
\begin{eqnarray*}
 \textrm{($\mbox{LR2}_\alpha$)}	&\min&
	\sum_{j\in J}\sum_{t\in {\cal T}}\sum_{s\in {\cal S}}\sum_{k\in K_{jt}^s}(\pi_s g_{jkt}^s-\alpha_{jt}^s)w_{jkt}^s+\sum_{j\in J}\sum_{t\in {\cal T}}\sum_{s\in {\cal S}}\sum_{k\in K_{jt}^{'s}}(\pi_s h_{jkt}^s+\alpha_{jt}^s)v_{jkt}^s \hspace{11ex} \\
	& & - \sum_{j\in J}\sum_{t\in {\cal T}}\sum_{s\in {\cal S}}\alpha_{jt}^sb_{jt}^s\\
	& & \\
	&\mbox{s.t.} &  \eqref{ww1}-\eqref{ww3},\eqref{wbound},\eqref{vbound}.
\end{eqnarray*}

The subproblem ($\mbox{LR1}_\alpha$) has the integrality property, as shown next.
As a consequence, the linear programming relaxation yields the optimal objective value of this subproblem.

\begin{prop}\label{flow}
	  Subproblem \textnormal{($\mbox{LR1}_\alpha$)} has the integrality property.
\end{prop}
\dem
To prove this result, we show that coefficient matrix of constraints for ($\mbox{LR1}_\alpha$), denoted by $A$, is totally unimodular (TU). 
The matrix form of this problem is given by
$$
\begin{array}{c}
B_1\begin{cases}
\end{cases}\\
B_2\begin{cases}
\\\\\\\\
\end{cases}\\
B_3\begin{cases}
\\\\\\\\
\end{cases}
\end{array}
\left(
\begin{array}{cccc:cccc}
{\bf T}&{\bf T}&\ldots&{\bf T}&{\bf -T^{\prime}}&{\bf -T^{\prime}}&\ldots&{\bf -T^{\prime}}\\
\hdashline
{\bf T}&{\bf \Theta}&\ldots&{\bf \Theta}&{\bf -T^{\prime}}&{\bf \Theta^{\prime}}&\ldots&{\bf \Theta^{\prime}}\\
{\bf \Theta}&{\bf T}&\ldots&{\bf \Theta}&{\bf \Theta^{\prime}}&{\bf -T^{\prime}}&\ldots&{\bf \Theta^{\prime}}\\
\vdots&\vdots&\ddots&\vdots&\vdots&\vdots&\ddots&\vdots\\
{\bf \Theta}&{\bf \Theta}&\ldots&{\bf T}&{\bf \Theta^{\prime}}&{\bf \Theta^{\prime}}&\ldots&{\bf -T^{\prime}}\\
\hdashline
-{\bf T}&{\bf \Theta}&\ldots&{\bf \Theta}&{\bf T^{\prime}}&{\bf \Theta^{\prime}}&\ldots&{\bf \Theta^{\prime}}\\
{\bf \Theta}&-{\bf T}&\ldots&{\bf \Theta}&{\bf \Theta^{\prime}}&{\bf T^{\prime}}&\ldots&{\bf \Theta^{\prime}}\\
\vdots&\vdots&\ddots&\vdots&\vdots&\vdots&\ddots&\vdots\\
{\bf \Theta}&{\bf \Theta}&\ldots&-{\bf T}&{\bf \Theta^{\prime}}&{\bf \Theta^{\prime}}&\ldots&{\bf T^{\prime}}\\
\end{array}
\right)
\begin{pmatrix}
{\bf z}_1\\
{\bf z}_2\\
\vdots\\
{\bf z}_m\\
\hdashline
{\bf z'}_1\\
{\bf z'}_2\\
\vdots\\
{\bf z'}_m\\
\end{pmatrix}
\leqslant
\begin{pmatrix}
{\bf p}\\
\hdashline
{\bf e^{\prime}}_1\\
{\bf e^{\prime}}_2\\
\vdots\\
{\bf e^{\prime}}_m\\
\hdashline
{\bf y'_{1}}\\
{\bf y'_{2}}\\
\vdots\\
{\bf y'_{m}}
\end{pmatrix}
$$
where
\begin{itemize}
\item ${\bf T}=
\begin{array}{c}
1\\
2\\
3\\
\vdots\\
T
\end{array}
\stackrel{
\begin{array}{ccccc}
1&2&3&\ldots&T\\
\end{array}
}{
\begin{pmatrix}
1&0&0&\ldots&0\\
1&1&0&\ldots&0\\
1&1&1&\ldots&0\\
\vdots&\vdots&\vdots&\ddots&\vdots\\
1&1&1&\ldots&1\\
\end{pmatrix}
}
$, 
$
{\bf T'}=
\begin{array}{c}
1\\
2\\
3\\
\vdots\\
T
\end{array}
\stackrel{
	\begin{array}{cccccccc}
	&1&&2&&3&\ldots&T-1\\
	\end{array}
}{
\begin{pmatrix}
0&0&0&\ldots&0\\
1&0&0&\ldots&0\\
1&1&0&\ldots&0\\
\vdots&\vdots&\vdots&\ddots&\vdots\\
1&1&1&\ldots&1\\
\end{pmatrix}
}
$.
\item
${\bf \Theta}$ and ${\bf \Theta}^{\prime}$ are  $T\times T$ and $T\times (T-1)$ null matrices, respectively.
\item 
${\bf z}_{i}=
\begin{pmatrix}
z_{i1}\\
\vdots\\
z_{iT}
\end{pmatrix}$
for $i\in I$,\hspace{0.4cm}
${\bf z}_{i}^{\prime}=
\begin{pmatrix}
z_{i1}^{\prime}\\
\vdots\\
z_{iT}^{\prime}
\end{pmatrix}$
for $i\in I$,\hspace{0.4cm}
${\bf p}=\begin{pmatrix}
	p_1-\sum_{i\in I} \bar{y}_{i0}\\
	\vdots\\
	p_T-\sum_{i\in I} \bar{y}_{i0}\\	
\end{pmatrix}$ integer.
\item ${\bf e}_i^{\prime}=
\begin{pmatrix}
e_i-\bar{y}_{i0}\\
\vdots\\
e_i-\bar{y}_{i0}
\end{pmatrix}
$
is a vector of $T$ integer components with identical value for $i\in I$.
\item 
$
{\bf y}_{i}^{\prime}=
\begin{pmatrix}
\bar{y}_{i0}\\
\vdots\\
\bar{y}_{i0}
\end{pmatrix}
$
is a vector of $T$ integer components with identical value for $i\in I$.
\end{itemize}
 
Observe that row-block $B_1$ is related with the family of constraints \eqref{ptrest}, $B_2$ with \eqref{ytub} and $B_3$ with \eqref{ytlb}.
In order to prove that $A$ is TU, applying Prop.\ 2.1 in \cite{Nemhauser}, we can remove 
the row-block $B_3$ because it may be obtained by duplicating the 
row-block $B_2$ and multiplying all its elements by $-1$.

Furthermore, 
${\bf T}=({\bf 1}|{\bf T'})$ where ${\bf 1}$ is a column vector of $T$ components with value $1$. 
Consequently, matrix
$$
 \left(
 \begin{array}{cccc}
 {\bf -T'}&{\bf -T'}&\ldots&{\bf -T'}\\
\hdashline
 {\bf -T'}&{\bf \Theta'}&\ldots&{\bf \Theta'}\\
 {\bf \Theta'}&{\bf -T'}&\ldots&{\bf \Theta'}\\
 \vdots&\vdots&\ddots&\vdots\\
 {\bf \Theta'}&{\bf \Theta'}&\ldots&{\bf -T'}
 \end{array}
 \right)
 $$
is the result of multiplying by $-1$ matrix
$$
  A'= (a'_{ij}) =
 \left(
 \begin{array}{cccc}
 {\bf T}&{\bf T}&\ldots&{\bf T}\\
\hdashline
 {\bf T}&{\bf \Theta}&\ldots&{\bf \Theta}\\
 {\bf \Theta}&{\bf T}&\ldots&{\bf \Theta}\\
 \vdots&\vdots&\ddots&\vdots\\
 {\bf \Theta}&{\bf \Theta}&\ldots&{\bf T}
 \end{array}
 \right)
 \begin{array}{l}
 \left.\right\}A'_0=\{ 1,\ldots ,T\}\\
 \left.\right\}A'_1=\{ T+1,\ldots ,2T\}\\
 \left.\right\}A'_2=\{ 2T+1,\ldots ,3T\}\\
 \vdots\\
 \left.\right\}A'_m=\{ mT+1,\ldots ,(m+1)T\}
 \end{array}
 $$
and removing columns $1,m+1,\ldots ,(T-1)m+1$. 

Therefore, 
in order to conclude the proof of this proposition it is only needed to prove that $A'$ is TU.

By \cite[Theorem 2.7]{Nemhauser},  $A'$ is TU if and only if for every subset $R$ of rows of $A'$ it is possible to find a partition $R=R_1 \cup R_2$ such that 
$$ \left|\sum_{i\in R_1}a'_{ij}-\sum_{i\in R_2}a'_{ij}\right|\leq 1, \quad \mbox{for each column $j$ of $A'$}.$$

Let $R=R'_0\cup\ldots\cup R'_m$ be a subset of rows of $A'$, 
such that $R'_i = \left\{k_{i1},\ldots,k_{i|R'_i|}\right\}\subseteq A'_i:=\{iT+1,\ldots,(i+1)T\}$
for $i\in\{0,\ldots,m\}$,  
with $k_{ij}\leq k_{i,j+1}$ for $j\in\{1,\ldots,|R_i^{\prime}|-1\}$. Let
\begin{eqnarray*}
R_i''& = & \{k_{i j} \in R_i' \: : \:  \mbox{with $j$ even}\},\\
R_i'''&=&  \{k_{i j} \in R_i' \: : \:  \mbox{with $j$ odd}\}.
\end{eqnarray*}
Note that it can be the case that $R_i'=R_i''=R_i'''=\emptyset$ for some values of $i$.

We have then that, 
for every $i=0,\ldots ,m$ and every column $j$ of $A'$,
vector $(a'_{\ell j})_{\ell \in R_i'}$ will have one of the shapes $(0,\ldots ,0)$, 
$(0,\ldots ,0,1,\ldots ,1)$ or $(1,\ldots,1)$. The first case is trivial, thus we will restrict out attention to the
second and third cases. Let $f_{ij}$ be the position of the first 1 in these types of vectors. If $f_{ij}$ is odd, the sum 
\begin{equation} \label{sumaalterna}
\sum_{\ell \in R'''_i} a'_{\ell j}-\sum_{\ell \in R_i''} a'_{\ell j}
\end{equation}
will take value 1 if the number of addends (given by $|R'_i|$) is odd or 0 if $|R'_i|$ is even. Similarly, 
if $f_{ij}$ is even, \eqref{sumaalterna} will take value -1 if $|R'_i|$ is even and 0 if $|R'_i|$ is odd.
Therefore
\begin{equation} \label{sumasalternas} 
\sum_{\ell \in R'''_i} a'_{\ell j}-\sum_{\ell \in R''_i} a'_{\ell j} \in 
\left\{ 
\begin{array}{cc}
\{ -1,0\} & \mbox{ if $|R'_i|$ is even},\\ 
\{ 0,1\}  & \mbox{ if $|R'_i|$ is odd}.
\end{array}
\right.
\end{equation}
Hence,  consider the following partition,
\begin{eqnarray*}
R_1&=& R_0''' \cup \bigcup_{i=1 \atop |R_0'|+|R_i'| \, \mbox{\tiny even}  }^m  R_i''\cup \bigcup_{i=1  \atop |R_0'|+|R_i'| \, \mbox{\tiny odd}  }^m  R_i''',\\  
R_2&=& R_0'' \cup \bigcup_{i=1, \atop |R_0'|+|R_i'| \, \mbox{\tiny even}  }^m  R_i'''\cup \bigcup_{i=1  \atop |R_0'|+|R_i'| \, \mbox{\tiny odd}  }^m  R_i''.\\  
\end{eqnarray*}
In order to get the value of the sum 
\begin{equation} \label{ultimasuma}
\sum_{\ell \in R_1} a'_{\ell j} -\sum_{\ell \in R_2} a'_{\ell j}
\end{equation}
for a given column $j$ of $A'$, we observe that all the 1s in a given column $j$ will correspond to only 
two row-blocks: $R'_0$ and any other $R'_s$. 
We distinguish four cases:
\begin{description}
\item[$|R'_0|$ is even and $|R's|$ is even.] 
In this case $R_1=R'''_0\cup R''_s$, $R_2=R''_0\cup R'''_s$ and \eqref{ultimasuma} can be calculated as
$$
(\sum_{\ell \in R'''_0} a'_{\ell j} - \sum_{\ell \in R''_0} a'_{\ell j}) - \\
(\sum_{\ell \in R'''_s} a'_{\ell j} - \sum_{\ell \in R''_s} a'_{\ell j}).
$$
Using \eqref{sumasalternas}, the contribution of both blocks will be 0 or -1. Hence, the difference will be then in the set $\{ -1,0,1\}$. 
\item[$|R'_0|$ is even and $|R's|$ is odd.] 
In this case $R_1=R'''_0\cup R'''_s$, $R_2=R''_0\cup R''_s$ and \eqref{ultimasuma} is
$$
(\sum_{\ell \in R'''_0} a'_{\ell j} - \sum_{\ell \in R''_0} a'_{\ell j}) + \\
(\sum_{\ell \in R'''_s} a'_{\ell j} - \sum_{\ell \in R''_s} a'_{\ell j}).
$$
Using \eqref{sumasalternas}, the contribution of 
block $R_0^{\prime}$ (respectively $R_s^{\prime}$)
 will be 
0 or -1 (resp.\  0 or 1). Hence, the sum will be then in the set $\{ -1,0,1\}$. 
\item[$|R'_0|$ is odd and $|R's|$ is even.] 
In this case $R_1=R'''_0\cup R'''_s$, $R_2=R''_0\cup R''_s$ and the reasoning of the previous case 
can be applied.
\item[$|R'_0|$ is odd and $|R's|$ is odd.] 
In this case $R_1=R'''_0\cup R''_s$, $R_2=R''_0\cup R'''_s$ and the situation is like in the first case.
\end{description}

Therefore, 
$$-1\le \sum_{\ell \in R_1} a'_{\ell j} -\sum_{\ell \in R_2} a'_{\ell j} \le 1,$$
and the result follows.
\fin

As far as subproblem ($\mbox{LR2}_\alpha$) is concerned, we first observe that it can be decomposed in such a way that a subproblem is obtained for each $j\in J$, $t\in {\cal T}$, $s\in {\cal S}$ as follows:
\begin{alignat}{3}
& \min & \quad & \sum_{k\in K_{jt}^s} (\pi_s g_{jkt}^s - \alpha_{jt}^s) w_{jkt}^s + \sum_{k\in K_{jt}^{\prime s}} (\pi_s h_{jkt}^s+\alpha_{jt}^s) v_{jkt}^s - \alpha_{jt}^s b_{jt}^s, & & \label{LR2-obj} 
\end{alignat}
\vspace{-7mm}
\begin{alignat}{3}
& \mbox{s.t.} & \quad & w_{j1t}^s+v_{j1t}^s\leq 1, & & \hspace{69mm} \label{LR2-constraint1} \\
& & & w_{jkt}^s\leq w_{j1t}^s, & &  k\in K_{jt}^s, \label{LR2-constraint2} \\
& & & v_{jkt}^s\leq v_{j1t}^s, & &  k\in K_{jt}^{'s}, \label{LR2-constraint3} \\
& & & w_{jkt}^s\in\{0,1\}, & & k\in K_{jt}^s, \label{LR2-constraint4} \\
& & & v_{jkt}^s\in\{0,1\}, & & k\in K_{jt}^{'s}. \label{LR2-constraint5}
\end{alignat}

In the following proposition we show that the integrality 
constraints \eqref{LR2-constraint4} and \eqref{LR2-constraint5} can be removed from the subproblem.

\begin{prop}
The matrix of coefficients of constraints  of subproblem  \eqref{LR2-obj}--\eqref{LR2-constraint5} is TU.
\end{prop}
\dem
For $j\in J$, $t\in {\cal T}$, $s\in {\cal S}$, linear constraints \eqref{LR2-constraint1}--\eqref{LR2-constraint3} can be represented as follows:
$$
\overbrace{
\left(
\begin{array}{ccccc:}
1&0&0&\ldots&0\\
-1&1&0&\ldots&0\\
-1&0&1&\ldots&0\\
\vdots&\vdots&\vdots&\ddots&\vdots\\
-1&0&0&\ldots&1\\
0&0&0&\ldots&0\\
0&0&0&\ldots&0\\
\vdots&\vdots&\vdots&\ddots&\vdots\\
0&0&0&\ldots&0\\
\end{array}
\right.
}^{C_1}
\overbrace{
	\left.
	\begin{array}{ccccc}
	 1&0&0&\ldots&0\\
	0&0&0&\ldots&0\\
	0&0&0&\ldots&0\\
	\vdots&\vdots&\vdots&\ddots&\vdots\\
	0&0&0&\ldots&0\\
    -1&1&0&\ldots&0\\
	-1&0&1&\ldots&0\\
	\vdots&\vdots&\vdots&\ddots&\vdots\\
	-1&0&0&\ldots&1\\
	\end{array}
	\right)
}^{C_2}
\begin{pmatrix}
w_{j1t}^s\\
w_{j2t}^s\\
w_{j3t}^s\\
\vdots \\
w_{j|K_{jt}^s|t}^s\\
\hdashline
v_{j1t}^s\\
v_{j2t}^s\\
v_{j3t}^s\\
\vdots\\
v_{j|K^{\prime s}_{jt}|t}^s\\
\end{pmatrix}
\leqslant
\begin{pmatrix}
1\\
0\\
0\\
\vdots \\
0\\
0\\
0\\
\vdots\\
0\\
\end{pmatrix}.
$$
We observe that each row of the coefficients matrix contains at most two nonzero elements. 
Besides, the columns of the matrix can be partitioned in 
two sets $C_1$ and $C_2$ such that two nonzero elements in a row are in the same set of columns if they have different signs, and in different sets of columns if they have the same sign. 
Since this property holds, by \cite[Corollary 2.8]{Nemhauser} the matrix is TU.
\fin

Besides being decomposable and having the unimodularity property, subproblem ($\mbox{LR2}_\alpha$) can be solved by inspection. 
To this end we note that, due to constraints \eqref{LR2-constraint1}--\eqref{LR2-constraint3}, in the optimal solution to any of the subproblems, $w$-variables and $v$-variables cannot take value one simultaneously. 
Instead, $w$-variables ($v$-variables) will take value one
if the sum of the negative coefficients associated with $w$-variables ($v$-variables) is the smallest one among both. 
For a fixed value of $j$, $t$ and $s$, we will take 
into account the following aspects:

\begin{enumerate}
	\item[(i)] For every $k\in K_{jt}^s$, 
	$\pi_s g_{jkt}^s - \alpha_{jt}^s < 0$ if $g_{jkt}^s < \frac{\alpha_{jt}^s}{\pi_s}$. 
	Since the last quotient does not depend on $j$ and parameters $g_{jkt}^s$ are sorted increasingly, if $\frac{\alpha_{jt}^s}{\pi_s}<g_{j|K_{jt}^s|t}^s$
	we will denote by $\ell$ the index such that	$\frac{\alpha_{jt}^s}{\pi_s} \in (g_{j \ell t}^s,g_{j,\ell+1,t}^s]$; 
	otherwise, we will take $\ell = |K_{jt}^s|$.
		
    Then, the binary values	of the $w$-variables that provide the lowest value of 
	$$\sum_{k\in K_{jt}^s}(\pi_s g_{jkt}^s-\alpha_{jt}^s)w_{jkt}^s$$ are: 
	$$w_{j1t}^s=\ldots=w_{j\ell t}^s=1,\;\; w_{j,\ell+1,t}^s=\ldots=w_{j|K_{jt}^s|t}^s=0, \quad \mbox{ if } 
	\frac{\alpha_{jt}^s}{\pi_s}<g_{j|K_{jt}^s|t}^s,
	$$ and 
	$$w_{j1t}^s=\ldots=w_{j |K_{jt}^s| t}^s=1, \quad \mbox{otherwise }.$$
	
	\item[(2)] 
	Analogously, the binary values of the $v$-variables that provide the lowest value of
	$$\sum_{k\in K_{jt}^{'s}}(\pi_s h_{jkt}^s+\alpha_{jt}^s)v_{jkt}^s$$ are
	$$v_{jkt}^s=\ldots=v_{j \ell^\prime t}^s=1,\;\;v_{j,\ell^\prime+1,t}^s=\ldots=v_{j|K_{jt}^{\prime s}|t}^s=0, \quad
	 \mbox{ if } 
	\frac{-\alpha_{jt}^s}{\pi_s}<h_{j|K_{jt}^{\prime s}|t}^s,$$ and 
	$$v_{jkt}^s=\ldots=v_{j |K_{jt}^{\prime s}| t}^s=1,\quad \mbox{otherwise,}$$
	where $\ell'$ denotes the index such that	$\frac{-\alpha_{jt}^s}{\pi_s} \in (h_{j \ell t}^s,h_{j,\ell+1,t}^s]$.
\end{enumerate}

Hence, each of the above subproblems can be solved straightforwardly as follows (we keep assuming $j$, $t$ and $s$ fixed in $J$, $T$ and $S$, respectively):

\medskip
We compare $\sum_{k\leq \ell}(\pi_s g_{jkt}^s-\alpha_{jt}^s)$ and 
$\sum_{k\leq \ell^\prime}(\pi_s h_{jkt}^s +\alpha_{jt}^s)$.
Then,
\begin{itemize}
	\item[ ] If 
	$\sum_{k\leq \ell}(\pi_s g_{jkt}^s-\alpha_{jt}^s) 
	\le \min \{ 0,\sum_{k\leq \ell^\prime}(\pi_s h_{jkt}^s +\alpha_{jt}^s)\}$ we set
	
	\smallskip
	\qquad $w_{j1t}^{*s}=\ldots=w_{j\ell t}^{*s}=1,\;\; w_{j,\ell+1,t}^{*s}=\ldots=w_{j|K_{jt}^s|t}^{*s}=0, \quad \mbox{ if } \ell < |K_{jt}^s|;$
	
	\qquad $w_{j1t}^{*s}=\ldots=w_{j |K_{jt}^s| t}^{*s}=1, \quad \mbox{ if } \ell = |K_{jt}^s|:$
	
	\qquad $v_{jkt}^{*s}=0,\;\; k\in K_{jt}^{'s}.$
	
	\item[ ] If  $\min \{ 0,\sum_{k\leq \ell}(\pi_s g_{jkt}^s-\alpha_{jt}^s)\} >\sum_{k\leq \ell'}(\pi_s h_{jkt}^s +\alpha_{jt}^s)$, we set
	
	\smallskip	
	\qquad 
	$v_{j1t}^{*s}=\ldots=v_{j \ell^\prime t}^{*s}=1
	,\;\;v_{j,\ell^\prime+1,t}^{*s}=\ldots=v_{j|K_{jt}^{\prime s}|t}^{*s}=0, \quad \mbox{ if } \ell^\prime < |K_{jt}^{\prime s}|;$
	
	\qquad $v_{j1t}^{*s}= \ldots=v_{j |K_{jt}^{\prime s}| t}^{*s}=1,\quad \mbox{ if } \ell^\prime = |K_{jt}^{\prime s}|;$

	\qquad $w_{jkt}^{*s}=0,\;\;k\in K_{jt}^{s}.$

	\item[ ] If  $\sum_{k\leq \ell}(\pi_s g_{jkt}^s-\alpha_{jt}^s)>0$ and $\sum_{k\leq \ell'}(\pi_s h_{jkt}^s +\alpha_{jt}^s)>0$,
	we set
	
	\qquad $w_{jkt}^{*s}=0,\;\;k\in K_{jt}^{s},$
	
	\qquad $v_{jkt}^{*s}=0,\;\;k\in K_{jt}^{'s}.$

\end{itemize}

Finally, a solution to ($\mbox{LR2}_\alpha$) is obtained also straightforwardly by considering all the solutions found above for all subproblems induced by $j \in J$, $t \in {\cal T}$ and $s \in {\cal S}$. 
The optimal value of subproblem ($\mbox{LR2}_\alpha$) is the sum of the optimal values of all such subproblems.

All the developments just presented show that the relaxed problem has the integrality property and thus the lower bound provided by our Lagrangian relaxation will not improve the lower bound provided by the linear relaxation of model (GMSCLP$^\prime$).
In other words, we are actually working with the Lagrangian relaxation of the linear relaxation of (GMSCLP$^\prime$).
However, we can use the solutions provided by subproblems ($\mbox{LR1}_\alpha$) and ($\mbox{LR2}_\alpha$) to obtain sharp upper bounds to model (GMSCLP$^\prime$).

\subsection{Deriving feasible solutions to (GMSCLP$^\prime$)}
\label{subsec:GSMC'-upperBounds}
Let us denote by $\{ \mathbf{z^*},\mathbf{z^{\prime *} \}}$ an integer optimal solution to (LR1$_\alpha$) for some $\alpha$. 
We note that such an integer solution can be easily obtained by solving (LR1$_{\alpha}$), for instance, using the Simplex method.
 
Indeed, since the coefficients matrix associated with 
{(LR1$_{\alpha}$)} is TU, this method guarantees that the obtained optimal solution 
is integer.

If we fix these variables in the original model (GMSCLP$^\prime$) the latter can be decoupled according to the different demand points, time periods and scenarios. 
Therefore, for every $j \in J$, $t \in {\cal T}$, and $s \in {\cal S}$ we can formulate the problem $\mbox{UB}(\mathbf{z^*},\mathbf{z^{\prime *}})_{jts}$ that consists of finding the optimal values of variables $w$ and $v$ as follows:

($\mbox{UB}(\mathbf{z^*},\mathbf{z^{\prime *}})_{jts}$)
\begin{alignat}{3}
& \min & \quad & \sum_{k\in K_{jt}^s}\pi_s g_{jkt}^sw_{jkt}^s+\sum_{k\in K_{jt}^{\prime s}}\pi_s h_{jkt}^s v_{jkt}^s & & \hspace{62mm} \nonumber 
\end{alignat}
\vspace{-8mm}
\begin{alignat}{3}
& \mbox{s.t.} & \quad & \sum_{k\in K_{jt}^s}w_{jkt}^s-\sum_{k\in K_{jt}^{\prime s}}v_{jkt}^s=M_{jts} & & \nonumber \\
& & & w_{j1t}^s+v_{j1t}^s\leq 1, & & \nonumber\\
& & & w_{jkt}^s\leq w_{j1t}^s, & &k\in K_{jt}^s \nonumber\\
& & & v_{jkt}^s\leq v_{j1t}^s, & &k\in K_{jt}^{\prime s},\nonumber\\
& & & w_{jkt}^s\in\{0,1\}, & & k\in K_{jt}^s,\nonumber\\
& & & v_{jkt}^s\in\{0,1\}, & & k\in K_{jt}^{\prime s}, \hspace{55mm} \nonumber
\end{alignat}
where $M_{jts}:=\sum_{i\in I}a_{ijt}^s\left( \bar{y}_{i0}+\sum_{\tau=1}^t z_{i\tau}^*-\sum_{\tau=1}^{t-1} z^{\prime *}_{i\tau}\right)-b_{jt}^s$.
It is important to remark that problem ($\mbox{UB}(\mathbf{z^*},\mathbf{z^{\prime *}})_{jts}$) is always feasible. 
In fact, taking \eqref{ptrest} into account as well as the binary nature of the parameters $a_{ijt}^s$ we know that
$$0 \le \sum_{i\in I}a_{ijt}^s\left( \bar{y}_{i0}+\sum_{\tau=1}^t z_{i\tau}^*-\sum_{\tau=1}^{t-1} z^{\prime *}_{i\tau}\right)\leq p_t,$$ 
thus, $M_{jts}$ is within the variation range of 
$\sum_{k\in K_{jt}^s}w_{jkt}^s-\sum_{k\in K_{jt}^{\prime s}}v_{jkt}^s$.
For instance, if we consider the extreme case
$M_{jts}= p_t-b_{jt}^s$, the solution will be $\widetilde{w}_{jkt}^{s}=1$ for $k \in K_{jt}^s$ and $\widetilde{v}_{jkt}^{s}=0$ for $k\in K_{jt}^{\prime s}$; 
in the other extreme we would have 
$M_{jts}=-b_{jt}^s$, 
and the corresponding solution would be
$\widetilde{v}_{jkt}^{s}=1$ for $k \in K_{jt}^{\prime s}$ and $\widetilde{w}_{jkt}^{s}=0$ for $k\in K_{jt}^{s}$.

In general, for every $j\in J$, $t\in {\cal T}$, $s\in {\cal S}$, 
\begin{itemize}
	\item[ ] if $M_{jts}<0$, the optimal solution to ($\mbox{UB}(\mathbf{z^*},\mathbf{z^{\prime *}})_{jts}$) is  
	
	\qquad $\widetilde{w}_{jkt}^{s} = 0, \quad k \in K_{jt}^s$;
	
	\qquad $\widetilde{v}_{jkt}^{s} = 1, \quad k \leq -M_{jts}$;
	
	\qquad $\widetilde{v}_{jkt}^{s} = 0, \quad k>-M_{jts}$;
	
	\item[ ] if $M_{jts}>0$, then the optimal solution becomes
	
	\qquad $\widetilde{v}_{jkt}^{s} = 0, \quad k\in K_{jt}^{'s}$;
	
	\qquad $\widetilde{w}_{jkt}^{s} = 1, \quad k\leq M_{jts}$;
	
	\qquad $\widetilde{w}_{jkt}^{s} = 0, \quad k>M_{jts}$;
	
	 \item[ ] if $M_{jts}=0$, then
	
	\qquad $\widetilde{w}_{jkt}^{s} = 0, \quad k \in K_{jt}^s$;
	
	\qquad $\widetilde{v}_{jkt}^{s} = 0, \quad k\in K_{jt}^{'s}$.

\end{itemize}

Hence, the optimal value of ($\mbox{UB}(\mathbf{z^*},\mathbf{z^{\prime *}})_{jts}$) is 
$\sum_{k=1}^{M_{jts}}\pi_s h_{jkt}^s$ if $M_{jts}<0$ and $\sum_{k=1}^{M_{jts}}\pi_s g_{jkt}^s$, otherwise. 

Combining the optimal solutions to the above subproblems we end up with an upper bound on the optimal value of model (GMSCLP$^\prime$) that can be computed in the following way:
\begin{multline}
\sum_{i\in I}\sum_{t\in {\cal T}}\left(o_{it}+\sum_{\tau=t}^{T}f_{i\tau}\right)z_{it}^*+
\sum_{i\in I}\sum_{t\in {\cal T}\setminus\{T\}}\left(c_{it}-\sum_{\tau=t+1}^{T}f_{i\tau}\right)z^{'*}_{it}\\
+
\sum_{i\in I}\sum_{t\in {\cal T}}f_{it} \bar{y}_{i0}+\sum_{j\in J}\sum_{t\in {\cal T}}\sum_{s\in {\cal S}} {\cal V}(\mbox{UB}(\mathbf{z^*},\mathbf{z^{\prime *}})_{jts}), \label{UB}
\end{multline}

\noindent
where ${\cal V}(\mbox{UB}(\mathbf{z^*},\mathbf{z^{\prime *}})_{jts})$ denotes the optimal value of subproblem $\mbox{UB}(\mathbf{z^*},\mathbf{z^{\prime *}})_{jts}$.

\subsection{Optimizing the Lagrangian dual}
\label{subsec:Lagrangian-dual}
We have now a mechanism for solving the  Lagrangian relaxed problem for some vector of multipliers as well as a procedure for deriving feasible solutions to the problem.

Our goal when using the Lagrangian relaxation is to optimize the Lagrangian dual, i.e., to solve $$\max_{\alpha} {\cal V}(\mbox{LR}_\alpha),$$
where ${\cal V}(\mbox{LR}_\alpha)$ {denotes} the optimal value of ($\mbox{LR}_\alpha$). 
This can be accomplished by using the subgradient method \citep[see, for instance,][]{Guignard2003}.
This is an iterative procedure such that in each iteration $\kappa$ we obtain a new vector of multipliers that we denote by $\alpha[\kappa]$, whose generic component is $\alpha_{jt}^s[\kappa]$, $j \in J$, $t \in {\cal T}$, $s \in {\cal S}$.

We initialize the process by setting $\alpha_{jt}^s[0]=0$, $j \in J$, $t \in {\cal T}$, $s \in {\cal S}$.
Then, we obtain  an initial lower bound, LB, by solving 
subproblems $(\mbox{LR1}_{\alpha[0]})$ and $(\mbox{LR2}_{\alpha[0]})$, and we obtain an initial upper bound, UB, by solving the corresponding  subproblems $\mbox{UB}(\mathbf{z^*},\mathbf{z^{\prime *}})_{jts}$.

In each generic iteration $\kappa$, subproblem $(\mbox{LR1}_{\alpha[\kappa]})$ is solved yielding an objective function value of ${\cal V}(\mbox{LR1}_{\alpha[\kappa]})$.
Subproblem $(\mbox{LR2}_{\alpha[\kappa]})$ is also solved rendering an objective function value of $\mbox{LR2}_{\alpha[\kappa]}$.
Hence, in each iteration $\kappa$, the lower bound produced is 
	$\mbox{LB}[\kappa] = {\cal V}(\mbox{LR1}_{\alpha[\kappa]}) + {\cal V}(\mbox{LR2}_{\alpha[\kappa]})$.
If $\mbox{LB}[\kappa] > \mbox{LB}$ then,  $\mbox{LB}$ is set equal to $\mbox{LB}[\kappa]$. 

To obtain an upper bound in an iteration $\kappa$, subproblems $\mbox{UB}(\mathbf{z^*},\mathbf{z^{\prime *}})_{jts}$ ($j \in J$, $t \in {\cal T}$, $s \in {\cal S}$) are solved using the values of $z$ and $z^\prime$ variables obtained from $\mbox{LR1}_{\alpha[\kappa]}$.
The solutions of these subproblems render a feasible solution to (GMSCLP$^\prime$) that can be evaluated using \eqref{UB}. 
If this value is smaller than the incumbent best upper bound then the latter is replaced by the former.

Finally, in each iteration, the Lagrangian multipliers are updated according to
$$\alpha_{jt}^s[\kappa+1] = \alpha_{jt}^s[\kappa] +  \varepsilon[\kappa] \: \cdot \: \frac{UB-\mbox{LB}[\kappa]}{||\gamma[\kappa]||^2} \cdot \gamma[\kappa],$$
where
$$\gamma[\kappa] = \sum_{i\in I}a_{ijt}^s\left( \bar{y}_{i0}+\sum_{\tau=1}^{t}z^{*}_{i\tau}-\sum_{\tau=1}^{t-1} z^{\prime*}_{i\tau}\right)-b_{jt}^s-\sum_{k\in K_{jt}^s}w_{jkt}^{*s}
+\sum_{k\in K_{jt}^{'s}}v_{jkt}^{*s}.$$
In  the expression of $\gamma[\kappa]$ an `*' next to a variable indicates the optimal value of the variable obtained in iteration $\kappa$.

It must be noted that $\varepsilon[\kappa]$ is a ``step'' parameter to be changed during the procedure.
In the numerical results section we provide the implementation details.

\section{Numerical experiments}
\label{sec:numerical_experiments}
In this section we report and discuss the numerical experiments conducted to evaluate the models and methodological developments proposed in the previous sections.

\subsection{Data generation}
\label{subsec:data_generation}
In order to perform the numerical experiments we coded an instance generator for obtaining the test instances.

We  generated instances assuming that the location of the customers define the potential location for the facilities. 
Accordingly, we assumed $|I|=|J|$. 
For these cardinalities we considered the values $5$, $10$, $30$, $50$ and $100$.
Additionally, we assumed $T=S=3$, that is, 3 periods in the planning horizon and 3 scenarios for capturing  the uncertainty.

For every value of  $|I|=|J|$ we have generated 5 instances according to the following methodology:

\begin{itemize}
	\item Each demand point is generated as a random point in the rectangle $[0,10]\times [0,50]$.  
	As mentioned above, such point 
	also defines a potential location for the facilities.
	
	\item For every $i\in I$ and $t\in {\cal T}$ the opening, closing and operating costs ($o_{it}$, $c_{it}$, $f_{it}$) 
	are randomly generated according to a 
	continuous uniform distribution in [1,10].
		
	\item For every $i\in I$, $j\in J$, $t\in {\cal T}$, $s\in {\cal S}$, the covering capability $a_{ijt}^s$ is obtained as follows: 
	
	\begin{itemize}
		\item First set a radius $r=8$. 
		Afterwards, for $i \in I$, $j \in J$,  $s\in {\cal S}$, we set $a_{ij1}^s=1$ if $d_{ij}\leq r$, where $d_{ij}$ denotes the  Euclidean distance between points $i$ and $j$.
		
		\item For the remaining periods we consider that the distance covered by a possible location is decreasing. 
		This is accomplished by decreasing $r$ by $20\%$ from one period to the following one. 
		
		In any case, for every $t=2,\ldots, T$ we set again $a_{ijt}^s=1$ for $s \in {\cal S}$ if $d_{ij}\leq r$. 
				
		\item Since we want the covering capabilities $a_{ijt}^s$ also to depend on the scenarios, the above values $a_{ijt}^s$ are then modified as follows:
		
		\begin{itemize}
			\item[ ] For every scenario $s \in {\cal S}$ we randomly select a number of locations given by round($0.2 \times m$) 
			and set  $a_{ijt}^s=0$, $i \in I_0$, $j \in J$, $t \in {\cal T}$.
			$I_0$ denotes the set of locations selected.
		\end{itemize}
	
	\end{itemize}
	
	\item For every $j\in J$, $t\in {\cal T}$, and $s\in {\cal S}$ we set 
	$b_{jt}^s=$round$(0.3\times \sum_{i\in I}a_{ijt}^s)$.
	
	\item For all $i\in I,$ we set $e_i=2$.
	
	\item For every $t \in {\cal T}$, $p_t$ is randomly generated according to a discrete uniform distribution $U\{\max \{1,0.1|J|\}, 0.3|J|\}$.
	
	\item 
	For fixed  values $j\in J$, $t\in T$ and $s \in {\cal S}$ 
	the costs $g^s_{jkt}$, $k \in\{1,\ldots,|K_{jt}^s|\}$, are generated as follows:
		we generate $|K_{jt}^s|$ values according to a continuous uniform distribution in the set $[-10,-1]$ and sort them non-decreasingly before assigning them to $g^s_{j1t},\dots,g^s_{j|K_{jt}^s|t}$, respectively.
		This way we ensure that $g_{jkt}^s \le g_{j,k+1,t}^s$ for $k=1,\ldots,|K_{jt}^s|-1$.
	
	\item
	For fixed  values $j\in J$, $t\in T$ and $s \in {\cal S}$, 
	the costs $h^s_{jkt}$, $k \in\{1,\ldots,|K_{jt}^{\prime s}\}$, are generated as follows:
		we generate $|K_{jt}^{\prime s}|$ values according to a continuous uniform distribution in the set $[1,10]$ and sort them non-decreasingly before assigning them to $h^s_{j1t},\dots,h^s_{j|K_{jt}^{\prime s}|t}$, respectively.
		This way we ensure that $h_{jkt}^s \le h_{j,k+1,t}^s$ for $k=1,\ldots,|K_{jt}^{\prime s}|-1$.
	
	\item Finally, concerning the probabilities $\pi_s$, $s \in {\cal S}$ we start by randomly generating them according to a continuous uniform distribution $U[0,1]$.
	Afterwards, 
	each such value is divided by their sum to ensure that $\sum_{s\in {\cal S}}\pi_s=1$.
\end{itemize}

We generated $5$ instances for each value of $|J|$, such that the set $J$ changes from one instance to another.

All the computational tests were performed using a machine 
Intel(R) Core(TM) i7-4790K CPU 32 GB RAM. Besides, models were coded in C++ using ILOG Concert Technology CPLEX 7.0.

We start by reporting the tests performed using the model proposed in Section~\ref{general_model}.
In Table~\ref{tab1},  the first column `LP-gap' is the percentage gap between the optimal value of the problem and the linear relaxation bound. Observe that the optimal value of (GMSCLP$^{\prime}$) could be negative or zero. In problems with possible negative or zero optimal values, zero cannot be used as a lower bound for the optimal value, which is implicit in the formulas usually considered for computing the gaps. Consequently, to 
calculate the relative gaps, an actual lower bound for (GMSCLP$^{\prime}$), that we will denote as LB$_{0}$, should be considered. One possibility for getting such bound is to solve the linear relaxation of (GMSCLP$^{\prime}$) removing the expected cost of shortage in the coverage from the formulation, i.e., we solve the linear relaxation of the following problem,
\begin{alignat}{3}
& \mbox{min} & \quad & \sum_{i\in I}\sum_{t\in {\cal T}}\left(o_{it}+\sum_{\tau=t}^{ T}f_{i\tau}\right)z_{it}+
\sum_{i\in I}\sum_{t\in {\cal T}\setminus\{ T\}}\left(c_{it}-\sum_{\tau=t+1}^{ T}f_{i\tau}\right)z^\prime_{it}+
\sum_{i\in I}\sum_{t\in {\cal T}}f_{it}{ \bar{y}_{i0}} & & \nonumber \\
& & &
+ \sum_{s\in {\cal S}}\pi_s\sum_{t\in {\cal T}}\sum_{j\in J}\sum_{k\in K_{jt}^s}g_{jkt}^s w_{jkt}^s.\nonumber\\
& \mbox{s.t.} &\quad &  \eqref{ww1}-\eqref{ww3},\eqref{zbound}-\eqref{vbound},\eqref{ptrest}-\eqref{ytlb}. & \qquad & \hspace{4cm}\nonumber 
\end{alignat}
The LP-gap is then calculated as $\frac{\mbox{OPT-LP}}{\mbox{OPT-LB}_{0}}\cdot 100$, where `OPT' is the optimal value of (GMSCLP$^{\prime}$) and `LP' is the linear relaxation objective value of (GMSCLP$^{\prime}$).

The column `LP-time' of Table \ref{tab1} reports the CPU time, in seconds, required to solve the linear relaxation and finally, `OPT-time' is the CPU time, again in seconds, required to solve model (GMSCLP$^{\prime}$).

A time limit of 3 hours was considered when solving model (GMSCLP$^\prime$). 
We mark with a `*' the running times corresponding to instances that were not solved in this time limit.

\begin{table}[htbp]
	\centering
	\caption{LP-gap, LP-time and OPT-time results using model (GMSCLP$^\prime$). }
	\renewcommand{\arraystretch}{0.95}
	\begin{adjustbox}{max width=1.00\textwidth}
	\begin{tabular}{rrrrr}
		\toprule
		&       & \multicolumn{3}{c}{T=3/S=3} \\
		m/n   &       & LP-gap & LP-time & OPT-time \\
		\cmidrule{1-1}\cmidrule{3-5}
		5     &       & 10.11 & $<$0.01 & 0.21 \\
		5     &       & 39.06 &$<$0.01  & 0.01 \\
		5     &       & 17.53 & $<$0.01 & 0.01 \\
		5     &       & 13.53 &$<$0.01  & 0.01 \\
		5     &       & 7.54  & $<$0.01 & $<$0.01 \\
		\hdashline
		10    &       & 8.51  & $<$0.01  & 0.01 \\
		10    &       & 28.37 & $<$0.01  & 0.21 \\
		10    &       & 20.73 &$<$0.01  & 0.03 \\
		10    &       & 31.54 & $<$0.01  & 0.02 \\
		10    &       & 20.49 & $<$0.01  & 0.01 \\
		\hdashline
		30    &       & 29.34 & 0.06  & 3.33 \\
		30    &       & 24.47 & 0.03  & 1.30 \\
		30    &       & 20.20 & 0.05  & 1.57 \\
		30    &       & 24.49 & 0.05  & 3.54 \\
		30    &       & 21.86 & 0.03  & 1.42 \\
		\hdashline
		50    &       & 22.42 & 0.17  & 1320.88 \\
		50    &       & 16.42 & 0.19  & 3.62 \\
		50    &       & 16.05 & 0.11  & 1.79 \\
		50    &       & 25.69 & 0.13  & 135.08 \\
		50    &       & 21.62 & 0.16  & 1716.52 \\
		\hdashline
		100   &       & 19.84 & 1.34  & 10811.70* \\
		100   &       & 26.29 & 1.11  & 10804.90* \\
		100   &       & 24.82 & 1.45  & 10809.40* \\
		100   &       & 23.73 & 1.31  & 10809.90* \\
		100   &       & 28.66 & 0.97  & 10809.40* \\
		\bottomrule
	\end{tabular}%
		\end{adjustbox}
\label{tab1}%
\end{table}%

\subsection{Lower and upper bounds produced by the Lagrangian relaxation}
\label{subsec:LR-bounds}
For the Lagrangian relaxation based approach presented in Section~\ref{sec:Lagrangian-heuristic} two variants were implemented. 
In both we take $\varepsilon[0]=1.5$.
Then, this parameter is updated when the objective value of the relaxation is not improved after some iterations. 
In particular, for variant 1, $\epsilon[\kappa]=\frac{\epsilon[\kappa-1]}{2}$ if the objective value has not been improved in $10$ iterations. 
In variant 2, $\epsilon[\kappa]$ is  updated in the same way, but when the number of iterations without improvement is $5$. For 
both variants, we check different stop criteria, namely: stop (i) after 50 iterations, (ii) after 150 iterations, (iii) after 500 iterations, and (iv) when $\epsilon[k]<0.005$. 
Besides, in all variants, we stop the Lagrangian heuristic if $\frac{\mbox{UB}-\mbox{LB}}{|\mbox{UB}|}100\leq 0.01$.

Tables~\ref{tab2}--\ref{tab5} report the results of the Lagrangian heuristic variants. 
In Table \ref{tab2} we can see the gap between  best upper $\mbox{(UB)}$ and lower $\mbox{(LB)}$ bounds obtained in the Lagrangian heuristic. 
These gaps have been computed as $\frac{\mbox{UB}-\mbox{LB}}{\mbox{UB-LB}_{0}}100$. 
The best results (boldfaced) are the ones corresponding to variants 1.iii and 2.iii.

Table~\ref{tab3} reports the percentage gap between the upper bound that we obtain using the Lagrangian heuristic and the optimal values obtained solving model (GMSCLP$^\prime$).
 
As observed, in all the variants, we obtain a gap smaller than $4\%$. 
Some of  the gaps obtained for the instances with $m=n=100$ are negative since the upper bounds are compared with the best solution obtained in $3$ hours. 
This indicates that the upper bound obtained using the Lagrangian relaxation based procedure is better than the value of the best feasible solution found by the solver in 3 hours.

The percentage gaps between the lower bound obtained when using the Lagrangian relaxation based approach and the linear programming relaxation of model (GMSCLP$^\prime$) are reported in Table \ref{tab4}. These gaps are calculated as $\frac{\mbox{LP-LB}}{\mbox{LP-LB}_{0}}100$.

It is remarkable that, in variant 1.iii, the solutions obtained using Lagrangian relaxation procedure are exactly the same as those provided by the linear programming relaxation. 
In general, all variants render good values. 

Finally, Table \ref{tab5} reports the running times of the Lagrangian relaxation based procedure. 
In this table we observe that the largest CPU time observed was 92 seconds. 
This is a clear indication that the Lagrangian relaxation allows us to obtain high-quality feasible solutions within a very small CPU time.

In summary, from all these tables, we can observe that the variants 1.iii, 1.iv, 2.iii and 2.iv are the ones that provide the best results in general. For this reason, in the next subsection, these variants have been additionally tested  for new instances increasing the number of scenarios and time periods.

\begin{table}[htbp]
	\centering
	\caption{Lagrangian relaxation based procedures: percentage gaps between the best lower and upper bounds found.} 
	\renewcommand{\arraystretch}{0.95}
	\begin{adjustbox}{max width=1.00\textwidth}
	\begin{tabular}{rrrrrrrrrrr}
			&   &  & & &  &  &  & &  \\
		\toprule 
		\multicolumn{11}{c}{$Gap_{LB/UB}$} \\
		m/n     &    & Variant 1.i & Variant 1.ii & Variant 1.iii & Variant 1.iv & &Variant 2.i & Variant 2.ii & Variant 2.iii & Variant 2.iv \\
		\cmidrule{1-1} \cmidrule{3-6} \cmidrule{8-11}
		5     &       & 11.48 & 10.12 & \textbf{10.11} & 10.12 &       & 10.25 & \textbf{10.11} & \textbf{10.11} & 10.13 \\
		&       & 40.88 & \textbf{39.06} & \textbf{39.06} & 39.14 &       & 39.30 & \textbf{39.06} & \textbf{39.06} & 39.13 \\
		&       & 19.37 & 17.55 & \textbf{17.53} & 17.55 &       & 17.74 & \textbf{17.53} & \textbf{17.53} & 17.55 \\
		&       & 16.01 & \textbf{13.53} & \textbf{13.53} & 13.56 &       & 14.25 & \textbf{13.53} & \textbf{13.53} & 13.56 \\
		&       & 7.71  & 7.54  & \textbf{7.54} & 7.56  &       & 7.64  & \textbf{7.54} & \textbf{7.54} & 7.55 \\
		\hdashline
		10    &       & 13.39 & 9.90  & \textbf{9.87} & 9.89  &       & 10.61 & 9.88  & \textbf{9.87} & 9.88 \\
		&       & 32.21 & 28.43 & \textbf{28.37} & 28.42 &       & 28.77 & \textbf{28.37} & \textbf{28.37} & 28.43 \\
		&       & 30.29 & 21.17 & \textbf{21.09} & 21.13 &       & 22.92 & \textbf{21.09} & \textbf{21.09} & 21.13 \\
		&       & 41.49 & 31.69 & \textbf{31.54} & 31.59 &       & 36.92 & \textbf{31.54} & \textbf{31.54} & 31.61 \\
		&       & 27.91 & 20.62 & \textbf{20.49} & 20.53 &       & 22.37 & 20.50 & \textbf{20.49} & 20.55 \\
		\hdashline
		30    &       & 36.71 & 30.04 & 29.91 & 29.98 &       & 31.13 & \textbf{29.85} & \textbf{29.85} & 29.91 \\
		&       & 38.84 & 27.58 & \textbf{26.50} & 26.57 &       & 30.51 & 27.02 & 27.00 & 27.07 \\
		&       & 24.96 & 21.22 & \textbf{21.05} & 21.09 &       & 22.80 & 21.06 & \textbf{21.05} & 21.09 \\
		&       & 35.55 & 25.77 & \textbf{25.46} & 25.51 &       & 27.54 & 25.48 & \textbf{25.46} & 25.53 \\
		&       & 28.76 & 23.13 & \textbf{22.92} & 22.97 &       & 24.88 & 22.93 & \textbf{22.92} & 22.98 \\
		\hdashline
		50    &       & 28.50 & 24.35 & \textbf{23.88} & 23.94 &       & 26.00 & 23.90 & \textbf{23.88} & 23.94 \\
		&       & 44.57 & 17.82 & \textbf{17.02} & 17.06 &       & 22.38 & 18.17 & 18.12 & 18.17 \\
		&       & 19.19 & 16.20 & \textbf{16.05} & 16.09 &       & 19.24 & 16.09 & \textbf{16.05} & 16.09 \\
		&       & 39.21 & 26.61 & \textbf{25.86} & 25.93 &       & 38.20 & 26.30 & 26.09 & 26.16 \\
		&       & 28.62 & 22.88 & \textbf{22.38} & 22.43 &       & 26.40 & 22.44 & \textbf{22.38} & 22.44 \\
		\hdashline
		100   &       & 26.78 & 22.64 & \textbf{22.51} & 22.57 &       & 23.85 & 22.53 & 22.52 & 22.57 \\
		&       & 37.78 & 26.51 & \textbf{26.26} & 26.32 &       & 27.53 & \textbf{26.26} & \textbf{26.26} & 26.32 \\
		&       & 34.22 & 23.66 & \textbf{23.32} & 23.38 &       & 29.40 & 23.38 & \textbf{23.32} & 23.38 \\
		&       & 28.96 & 23.81 & \textbf{23.53} & 23.60 &       & 26.52 & 23.84 & 23.81 & 23.87 \\
		&       & 40.05 & 28.03 & \textbf{27.72} & 27.80 &       & 30.59 & 27.77 & 27.73 & 27.81 \\
		\bottomrule
	\end{tabular}%
		  \end{adjustbox}
\label{tab2}%
\end{table}%

\begin{table}[htbp]
	\centering
	\caption{Lagrangian relaxation based procedure: percentage gaps of the upper bound w.r.t. the optimal value of the problem.} 
	\renewcommand{\arraystretch}{0.95}
	\begin{adjustbox}{max width=1.00\textwidth}
	\begin{tabular}{rrrrrrrrrrr}
			&&&&&&&&&&\\
		\toprule
		\multicolumn{11}{c}{$Gap_{UB/OPT}$} \\
		\midrule
		m/n  &         & Variant 1.i & Variant 1.ii & Variant 1.iii & Variant 1.iv & &Variant 2.i & Variant 2.ii & Variant 2.iii & Variant 2.iv \\
		\cmidrule{1-1} \cmidrule{3-6} \cmidrule{8-11}
		5     &       & \textbf{$<$0.01} & \textbf{$<$0.01} & \textbf{$<$0.01} & \textbf{$<$0.01} &       & \textbf{$<$0.01} & \textbf{$<$0.01} & \textbf{$<$0.01} & \textbf{$<$0.01} \\
		&       & \textbf{$<$0.01} & \textbf{$<$0.01} & \textbf{$<$0.01} & \textbf{$<$0.01} &       & \textbf{$<$0.01} & \textbf{$<$0.01} & \textbf{$<$0.01} & \textbf{$<$0.01} \\
		&       & \textbf{$<$0.01} & \textbf{$<$0.01} & \textbf{$<$0.01} & \textbf{$<$0.01} &       & \textbf{$<$0.01} & \textbf{$<$0.01} & \textbf{$<$0.01} & \textbf{$<$0.01} \\
		&       & \textbf{$<$0.01} & \textbf{$<$0.01} & \textbf{$<$0.01} & \textbf{$<$0.01} &       & \textbf{$<$0.01} & \textbf{$<$0.01} & \textbf{$<$0.01} & \textbf{$<$0.01} \\
		&       & \textbf{$<$0.01} & \textbf{$<$0.01} & \textbf{$<$0.01} & \textbf{$<$0.01} &       & \textbf{$<$0.01} & \textbf{$<$0.01} & \textbf{$<$0.01} & \textbf{$<$0.01} \\
		\hdashline
		10    &       & \textbf{1.49} & \textbf{1.49} & \textbf{1.49} & \textbf{1.49} &       & \textbf{1.49} & \textbf{1.49} & \textbf{1.49} & \textbf{1.49} \\
		&       & \textbf{$<$0.01} & \textbf{$<$0.01} & \textbf{$<$0.01} & \textbf{$<$0.01} &       & \textbf{$<$0.01} & \textbf{$<$0.01} & \textbf{$<$0.01} & \textbf{$<$0.01} \\
		&       & \textbf{0.46} & \textbf{0.46} & \textbf{0.46} & \textbf{0.46} &       & \textbf{0.46} & \textbf{0.46} & \textbf{0.46} & \textbf{0.46} \\
		&       & \textbf{$<$0.01} & \textbf{$<$0.01} & \textbf{$<$0.01} & \textbf{$<$0.01} &       & \textbf{$<$0.01} & \textbf{$<$0.01} & \textbf{$<$0.01} & \textbf{$<$0.01} \\
		&       & \textbf{$<$0.01} & \textbf{$<$0.01} & \textbf{$<$0.01} & \textbf{$<$0.01} &       & \textbf{$<$0.01} & \textbf{$<$0.01} & \textbf{$<$0.01} & \textbf{$<$0.01} \\
		\hdashline
		30    &       & 0.80  & 0.80  & 0.80  & \textbf{0.71} &       & \textbf{0.71} & \textbf{0.71} & \textbf{0.71} & \textbf{0.71} \\
		&       & \textbf{2.68} & \textbf{2.68} & \textbf{2.68} & 3.35  &       & 3.35  & 3.35  & 3.35  & 3.35 \\
		&       & 1.50  & \textbf{1.07} & \textbf{1.07} & \textbf{1.07} &       & \textbf{1.07} & \textbf{1.07} & \textbf{1.07} & \textbf{1.07} \\
		&       & \textbf{1.29} & \textbf{1.29} & \textbf{1.29} & \textbf{1.29} &       & \textbf{1.29} & \textbf{1.29} & \textbf{1.29} & \textbf{1.29} \\
		&       & \textbf{1.36} & \textbf{1.36} & \textbf{1.36} & \textbf{1.36} &       & \textbf{1.36} & \textbf{1.36} & \textbf{1.36} & \textbf{1.36} \\
		\hdashline
		50    &       & \textbf{1.89} & \textbf{1.89} & \textbf{1.89} & \textbf{1.89} &       & \textbf{1.89} & \textbf{1.89} & \textbf{1.89} & \textbf{1.89} \\
		&       & 5.88  & \textbf{0.70} & \textbf{0.70} & 1.98  &       & 2.44  & 1.98  & 1.98  & 1.98 \\
		&       & \textbf{$<$0.01} & \textbf{$<$0.01} & \textbf{$<$0.01} & \textbf{$<$0.01} &       & \textbf{$<$0.01} & \textbf{$<$0.01} & \textbf{$<$0.01} & \textbf{$<$0.01} \\
		&       & 0.48  & \textbf{0.23} & \textbf{0.23} & 0.47  &       & 0.52  & 0.47  & 0.47  & 0.47 \\
		&       & 2.93  & \textbf{0.97} & \textbf{0.97} & \textbf{0.97} &       & 1.25  & \textbf{0.97} & \textbf{0.97} & \textbf{0.97} \\
		\hdashline
		100   &       & 3.74  & \textbf{3.34} & \textbf{3.34} & \textbf{3.34} &       & 3.51  & \textbf{3.34} & \textbf{3.34} & \textbf{3.34} \\
		&       & 0.40  & \textbf{-0.05} & \textbf{-0.05} & \textbf{-0.05} &       & \textbf{-0.05} & \textbf{-0.05} & \textbf{-0.05} & \textbf{-0.05} \\
		&       & \textbf{-2.00} & \textbf{-2.00} & \textbf{-2.00} & \textbf{-2.00} &       & \textbf{-2.00} & \textbf{-2.00} & \textbf{-2.00} & \textbf{-2.00} \\
		&       & 0.71  & \textbf{-0.27} & \textbf{-0.27} & 0.09  &       & 0.58  & 0.09  & 0.09  & 0.09 \\
		&       & \textbf{-1.32} & \textbf{-1.32} & \textbf{-1.32} & \textbf{-1.32} &       & \textbf{-1.32} & \textbf{-1.32} & \textbf{-1.32} & \textbf{-1.32} \\
		\bottomrule
	\end{tabular}%
\end{adjustbox}
\label{tab3}%
\end{table}%

\begin{table}[htbp]
	\centering
	\centering
	\caption{Lagrangian relaxation  versus linear programming relaxation: comparison of lower bounds.}
	\renewcommand{\arraystretch}{0.95}
	\begin{adjustbox}{max width=1.00\textwidth}
	\begin{tabular}{rrrrrrrrrrr}
		&&&&&&&&&&\\
		\toprule
		\multicolumn{11}{c}{$Gap_{LP/{  LB}}$} \\
		m/n      &   & Variant 1.i & Variant 1.ii & Variant 1.iii & Variant 1.iv & &Variant 2.i & Variant 2.ii & Variant 2.iii & Variant 2.iv \\
		\cmidrule{1-1} \cmidrule{3-6} \cmidrule{8-11}
		5     &       & 1.52  & 0.01  & \textbf{$<$0.01} & 0.01  &       & 0.16  & \textbf{$<$0.01} & \textbf{$<$0.01} & 0.02 \\
		&       & 2.99  & 0.01  & \textbf{$<$0.01} & 0.14  &       & 0.40  & \textbf{$<$0.01} & \textbf{$<$0.01} & 0.12 \\
		&       & 2.23  & 0.02  & \textbf{$<$0.01} & 0.02  &       & 0.26  & \textbf{$<$0.01} & \textbf{$<$0.01} & 0.02 \\
		&       & 2.87  & \textbf{$<$0.01} & \textbf{$<$0.01} & 0.04  &       & 0.84  & \textbf{$<$0.01} & \textbf{$<$0.01} & 0.04 \\
		&       & 0.19  & \textbf{$<$0.01} & \textbf{$<$0.01} & 0.02  &       & 0.11  & \textbf{$<$0.01} & \textbf{$<$0.01} & 0.02 \\
			\hdashline
		10    &       & 3.90  & 0.03  & \textbf{$<$0.01} & 0.02  &       & 0.81  & 0.01  & \textbf{$<$0.01} & 0.01 \\
		&       & 5.36  & 0.09  & \textbf{$<$0.01} & 0.08  &       & 0.56  & \textbf{$<$0.01} & \textbf{$<$0.01} & 0.08 \\
		&       & 11.66 & 0.10  & \textbf{$<$0.01} & 0.05  &       & 2.32  & 0.01  & \textbf{$<$0.01} & 0.05 \\
		&       & 14.54 & 0.22  & \textbf{$<$0.01} & 0.08  &       & 7.87  & 0.01  & \textbf{$<$0.01} & 0.10 \\
		&       & 9.32  & 0.16  & \textbf{$<$0.01} & 0.04  &       & 2.36  & \textbf{$<$0.01} & \textbf{$<$0.01} & 0.07 \\
			\hdashline
		30    &       & 9.70  & 0.18  & \textbf{$<$0.01} & 0.10  &       & 1.83  & 0.01  & \textbf{$<$0.01} & 0.09 \\
		&       & 16.79 & 1.47  & \textbf{$<$0.01} & 0.10  &       & 4.80  & 0.02  & \textbf{$<$0.01} & 0.09 \\
		&       & 4.54  & 0.22  & \textbf{$<$0.01} & 0.05  &       & 2.22  & 0.02  & \textbf{$<$0.01} & 0.06 \\
		&       & 13.54 & 0.42  & \textbf{$<$0.01} & 0.07  &       & 2.79  & 0.03  & \textbf{$<$0.01} & 0.09 \\
		&       & 7.57  & 0.27  & \textbf{$<$0.01} & 0.06  &       & 2.54  & 0.01  & \textbf{$<$0.01} & 0.07 \\
			\hdashline
		50    &       & 6.06  & 0.61  & \textbf{$<$0.01} & 0.07  &       & 2.78  & 0.02  & \textbf{$<$0.01} & 0.07 \\
		&       & 29.54 & 0.97  & \textbf{0.01} & 0.06  &       & 4.80  & 0.11  & 0.05  & 0.11 \\
		&       & 3.74  & 0.18  & \textbf{$<$0.01} & 0.05  &       & 3.80  & 0.05  & \textbf{$<$0.01} & 0.05 \\
		&       & 17.80 & 1.02  & \textbf{$<$0.01} & 0.09  &       & 16.41 & 0.35  & 0.06  & 0.16 \\
		&       & 6.17  & 0.64  & \textbf{$<$0.01} & 0.06  &       & 4.91  & 0.08  & $<$0.01  & 0.08 \\
			\hdashline
		100   &       & 5.11  & 0.16  & \textbf{$<$0.01} & 0.08  &       & 1.55  & 0.03  & 0.01  & 0.08 \\
		&       & 15.25 & 0.35  & \textbf{$<$0.01} & 0.09  &       & 1.73  & 0.01  & \textbf{$<$0.01} & 0.09 \\
		&       & 14.22 & 0.45  & \textbf{$<$0.01} & 0.08  &       & 7.93  & 0.08  & \textbf{$<$0.01} & 0.08 \\
		&       & 6.19  & 0.37  & \textbf{0.01} & 0.09  &       & 3.10  & 0.06  & 0.02  & 0.10 \\
		&       & 17.06 & 0.43  & \textbf{$<$0.01} & 0.11  &       & 3.97  & 0.08  & 0.02  & 0.12 \\
		\bottomrule
	\end{tabular}%
		\end{adjustbox}
\label{tab4}%
\end{table}%


\begin{table}[tbp]
	\centering
	\caption{ Lagrangian relaxation based procedure: CPU time (seconds).} 
	\renewcommand{\arraystretch}{0.95}
	\begin{adjustbox}{max width=1.00\textwidth}
\begin{tabular}{rrrrrrrrrrr}
	&&&&&&&&&&\\
	\toprule
	\multicolumn{11}{c}{ Lagrangian relaxation times} \\
	m/n   &   & Variant 1.i & Variant 1.ii & Variant 1.iii & Variant 1.iv & &Variant 2.i & Variant 2.ii & Variant 2.iii & Variant 2.iv \\
	\cmidrule{1-1} \cmidrule{3-6} \cmidrule{8-11}
	5     &       & 0.30  & 1.20  & 4.37  & 0.94  &       & \textbf{0.26} & 0.84  & 2.28  & 1.10 \\
&       & \textbf{0.28} & 0.86  & 1.77  & 0.57  &       & 0.30  & 0.62  & 2.03  & 0.34 \\
&       & 0.39  & 0.63  & 1.64  & 0.54  &       & \textbf{0.27} & 0.58  & 2.03  & 0.36 \\
&       & 0.27  & 0.94  & 1.65  & 0.47  &       & \textbf{0.25} & 0.50  & 2.20  & 0.47 \\
&       & \textbf{0.23} & 1.67  & 1.63  & 0.44  &       & 0.26  & 0.52  & 2.06  & 0.46 \\
	\hdashline
10    &       & \textbf{0.30} & 0.83  & 2.40  & 0.90  &       & \textbf{0.30} & 0.72  & 2.78  & 0.71 \\
&       & \textbf{0.30} & 0.92  & 2.44  & 0.99  &       & 0.32  & 0.73  & 2.82  & 0.46 \\
&       & 0.53  & 0.96  & 2.71  & 1.07  &       & \textbf{0.29} & 0.73  & 2.80  & 0.67 \\
&       & 0.29  & 0.76  & 2.50  & 0.97  &       & \textbf{0.25} & 0.75  & 2.53  & 0.55 \\
&       & 0.29  & 0.72  & 2.47  & 0.98  &       & \textbf{0.25} & 0.71  & 2.43  & 0.67 \\
	\hdashline
30    &       & \textbf{1.17} & 3.38  & 11.25 & 4.06  &       & \textbf{1.17} & 3.42  & 10.98 & 2.40 \\
&       & 1.22  & 3.23  & 10.97 & 5.12  &       & \textbf{1.11} & 3.37  & 10.94 & 2.71 \\
&       & 1.21  & 3.33  & 11.56 & 4.99  &       & \textbf{1.14} & 3.35  & 10.90 & 2.90 \\
&       & \textbf{1.12} & 3.33  & 10.93 & 4.46  &       & \textbf{1.12} & 3.39  & 10.91 & 2.68 \\
&       & \textbf{1.07} & 3.19  & 10.86 & 4.13  &       & 1.10  & 3.30  & 10.77 & 2.57 \\
	\hdashline
50    &       & 2.54  & 7.43  & 25.00 & 10.09 &       & \textbf{2.51} & 7.47  & 24.59 & 5.84 \\
&       & 2.61  & 7.80  & 25.69 & 14.65 &       & \textbf{2.56} & 7.70  & 24.98 & 7.68 \\
&       & \textbf{2.53} & 7.66  & 25.03 & 10.22 &       & \textbf{2.53} & 7.62  & 24.66 & 7.57 \\
&       & \textbf{2.51} & 7.63  & 25.36 & 12.07 &       & 2.54  & 7.65  & 25.06 & 8.41 \\
&       & \textbf{2.50} & 7.45  & 25.17 & 11.08 &       & \textbf{2.50} & 7.47  & 24.59 & 7.48 \\
	\hdashline
100   &       & \textbf{9.33} & 27.97 & 92.73 & 33.85 &       & 9.44  & 27.93 & 91.83 & 22.35 \\
&       & \textbf{8.98} & 26.78 & 88.44 & 33.70 &       & 9.04  & 26.57 & 87.22 & 19.33 \\
&       & \textbf{9.06} & 27.21 & 90.34 & 37.25 &       & 9.91  & 27.39 & 89.03 & 27.20 \\
&       & \textbf{9.06} & 27.04 & 89.47 & 37.46 &       & 9.24  & 27.19 & 88.20 & 24.85 \\
&       & 9.51  & 27.12 & 89.36 & 37.97 &       & \textbf{9.15} & 26.93 & 88.08 & 24.87 \\
	\bottomrule
\end{tabular}%
\end{adjustbox}
	\label{tab5}%
\end{table}%

\subsection{Lagrangian relaxation based procedure: additional tests}
\label{subsec:intensive-testing}
The very good results obtained so far for the Lagrangian relaxation based procedure, encouraged us to deepen the numerical tests by the inclusion of larger instances namely in terms of the number of scenarios and the number of time periods. 
With this purpose, we generated  an additional subset of instances using the following combinations of parameters $(T,S)=(3,5)$, $(T,S)=(5,5)$, $(T,S)=(5,10)$ and $(T,S)=(10,10)$ and using $m=n=5,10,30,50,100$.
The methodology for generating the instances is the same explained in Section~\ref{subsec:data_generation}.

Table \ref{tablbub1} reports the percentage gaps between the best upper bounds and lower bounds obtained when using the Lagrangian relaxation based procedure.  As before, these gaps are calculated as  $\frac{\mbox{UB}-\mbox{LB}}{\mbox{UB-LB}_{0}}100$.
Like in the previous tables, the values in bold indicate the minimum gap across all the studied variants. 
Looking into these tables we conclude that the best gaps are provided by variant 1.iii. 
We note that, in this variant, 500 iterations are executed in total updating $\epsilon$ after 10 iterations without improvement.
 
In Table~\ref{tabubopt1} we report the gap between the best upper bounds obtained by our Lagrangian heuristic and the optimal solutions obtained by solving model  (GMSCLP$^\prime$).These gaps are again calculated as $\frac{\mbox{UB-OPT}}{\mbox{UB-LB}_{0}}100$.

Regarding these gaps, all the variants exhibit a similar behavior, providing a small gap with respect to the optimal solution. 
In particular, the gaps obtained for $(T,S)=(3,5)$ are all smaller than $4\%$; the gaps for $(T,S)=(5,5)$ are smaller than $3\%$; the gaps for $(T,S)=(5,10)$ are smaller than $ 5\%$; and the gaps for $(T,S)=(10,10)$ are smaller than $5\%$. 
Consequently, we can say that the feasible solutions obtained using the Lagrangian heuristic are still good for larger values of $T$ and $S$. 
It must be noticed that, as before, the negative gaps are due to the cases where the  solver could not get in 3 hours a solution to model (GMSCLP$^\prime$).
In these cases, the comparison is made between the best feasible solution obtained using the Lagrangian heuristic and the best feasible solution found by the solver (in 3 hours).
 
With respect to the gaps between the lower bounds produced by Lagrangian relaxation and those provided by the linear programming relaxation, we can see in Table~\ref{tablplag1} that they are very small in all cases. Furthermore, looking into this table, the maximum gap is smaller than $ 0.2\%$. 
Although all variants present good results in terms of these gaps, we can emphasize variant 1.iii that $0\%$ gap in most of the cases.
 
{In addtion}, Tables~\ref{tabtimes1} and \ref{tabtimes2} report the CPU time in seconds required by the Lagrangian relaxation based procedure for the large instances tested. 
In these tables we also present the CPU time required by the solver for tackling model (GMSCLP$^\prime$).

For $(T,S)=(3,5)$ the largest 
running times of variant 2.iv are smaller than 54 seconds. 
The 
largest times of this same variant and $(T,S)=(5,5)$ are approximately 113 seconds. 
In the  case of $(T,S)=(5,10)$, the 
largest times are 250 seconds and for $(T,S)=(10,10)$, the 
largest running times are 815 seconds. 
It is remarkable that these CPU times are small if we compare them with the CPU time required by the solver for handling the model. 
Hence, we may conclude that the Lagrangian heuristic provides high-quality feasible solutions  in rather small CPU times for these values of $T$ and $S$ and $m=n\geq 30$.

\begin{sidewaystable}[htbp]
	\centering
	\caption{Lagrangian relaxation based procedures: percentage gaps between the best lower and upper bounds found for $T=3/S=5$, $T=S=5$, $T=5/S=10$ and $T=S=10$.}
	\renewcommand{\arraystretch}{1}
	\begin{adjustbox}{max width=1.10\textwidth}
		\begin{tabular}{rrrrrrrrrrrrrrrrrrrrr}
			\toprule
			\multicolumn{21}{c}{Gap$_{LB/UB}$} \\
			&       & \multicolumn{4}{c}{T=3/S=5}   &       & \multicolumn{4}{c}{T=5/S=5}   &       & \multicolumn{4}{c}{T=5/S=10}  &       & \multicolumn{4}{c}{T=10/S=10} \\
			m/n   &       & Variant 1.iii & Variant 1.iv & Variant 2.iii & Variant 2.iv &       & Variant 1.iii & Variant 1.iv & Variant 2.iii & Variant 2.iv &       & Variant 1.iii & Variant 1.iv & Variant 2.iii & Variant 2.iv &       & Variant 1.iii & Variant 1.iv & Variant 2.iii & Variant 2.iv \\
			\cmidrule{1-1}	\cmidrule{3-6}	\cmidrule{8-11}	\cmidrule{13-16} 	\cmidrule{18-21}
			5     &       & \textbf{15.60} & 15.62 & \textbf{15.60} & 15.62 &       & \textbf{18.48} & 18.50 & \textbf{18.48} & 18.52 &       & \textbf{4.56} & 4.58  & \textbf{4.56} & 4.58  &       & \textbf{32.73} & 32.78 & \textbf{32.73} & 32.81 \\
			&       & \textbf{6.42} & 6.43  & \textbf{6.42} & 6.43  &       & \textbf{18.66} & 18.69 & \textbf{18.66} & 18.68 &       & \textbf{11.35} & 11.38 & \textbf{11.35} & 11.38 &       & \textbf{15.42} & 15.44 & \textbf{15.42} & 15.45 \\
			&       & \textbf{9.50} & 9.51  & \textbf{9.50} & 9.51  &       & \textbf{18.84} & 18.85 & \textbf{18.84} & 18.85 &       & \textbf{4.28} & \textbf{4.28} & \textbf{4.28} & \textbf{4.28} &       & \textbf{6.02} & 6.04  & \textbf{6.02} & 6.04 \\
			&       & \textbf{0.23} & \textbf{0.23} & \textbf{0.23} & \textbf{0.23} &       & \textbf{0.01} & \textbf{0.01} & \textbf{0.01} & \textbf{0.01} &       & \textbf{22.80} & 22.83 & \textbf{22.80} & 22.83 &       & \textbf{18.93} & 18.96 & \textbf{18.93} & 18.97 \\
			&       & \textbf{6.69} & 6.70  & \textbf{6.69} & 6.71  &       & \textbf{0.01} & \textbf{0.01} & \textbf{0.01} & \textbf{0.01} &       & \textbf{13.32} & 13.33 & \textbf{13.32} & 13.35 &       & \textbf{12.97} & 12.98 & \textbf{12.97} & 12.99 \\
				\hdashline
			10    &       & \textbf{24.10} & 24.14 & \textbf{24.10} & 24.15 &       & \textbf{6.58} & 6.59  & \textbf{6.58} & 6.59  &       & \textbf{15.97} & 16.00 & \textbf{15.97} & 16.00 &       & \textbf{20.38} & 20.41 & 18.53 & 18.56 \\
			&       & \textbf{20.13} & 20.15 & \textbf{20.13} & 20.16 &       & \textbf{20.19} & 20.22 & \textbf{20.19} & 20.24 &       & \textbf{17.89} & 17.91 & \textbf{17.89} & 17.92 &       & \textbf{14.48} & 14.50 & \textbf{14.48} & 14.50 \\
			&       & \textbf{11.67} & 11.69 & \textbf{11.67} & 11.69 &       & \textbf{15.61} & 15.64 & \textbf{15.61} & 15.64 &       & \textbf{17.88} & 17.91 & \textbf{17.88} & 17.92 &       & \textbf{19.50} & 19.53 & \textbf{19.50} & 19.55 \\
			&       & \textbf{8.20} & 8.21  & \textbf{8.20} & 8.22  &       & \textbf{12.17} & 12.19 & \textbf{12.17} & 12.20 &       & \textbf{16.34} & 16.37 & \textbf{16.34} & 16.36 &       & \textbf{25.56} & 25.61 & \textbf{25.56} & 25.60 \\
			&       & \textbf{35.71} & 35.75 & \textbf{35.71} & 35.76 &       & \textbf{20.65} & 20.68 & \textbf{20.65} & 20.68 &       & \textbf{28.48} & 28.53 & \textbf{28.48} & 28.54 &       & \textbf{25.91} & 25.95 & \textbf{25.91} & 25.97 \\
				\hdashline
			30    &       & \textbf{24.90} & 24.96 & \textbf{24.90} & 24.97 &       & \textbf{19.01} & 19.05 & 19.48 & 19.53 &       & \textbf{18.29} & 18.33 & \textbf{18.29} & 18.33 &       & \textbf{16.97} & 17.01 & 17.04 & 17.08 \\
			&       & \textbf{21.39} & 21.45 & \textbf{21.39} & 21.45 &       & \textbf{23.23} & 23.27 & \textbf{23.23} & 23.29 &       & \textbf{25.25} & 25.30 & \textbf{25.25} & 25.31 &       & \textbf{19.44} & 19.49 & \textbf{19.44} & 19.80 \\
			&       & \textbf{21.11} & 21.16 & \textbf{21.11} & 21.16 &       & \textbf{21.67} & 21.73 & 21.69 & 21.75 &       & \textbf{19.99} & 20.04 & \textbf{19.99} & 20.04 &       & \textbf{20.86} & 20.91 & \textbf{20.86} & 20.91 \\
			&       & \textbf{24.71} & \textbf{24.77} & \textbf{24.71} & \textbf{24.77} &       & \textbf{17.13} & 17.17 & 17.15 & 17.20 &       & \textbf{17.63} & 17.67 & \textbf{17.63} & 17.67 &       & \textbf{16.09} & 16.13 & 16.11 & 16.15 \\
			&       & \textbf{24.44} & 24.50 & \textbf{24.44} & 24.51 &       & \textbf{12.37} & 12.40 & \textbf{12.37} & 12.40 &       & \textbf{19.77} & 19.82 & \textbf{19.77} & 19.82 &       & \textbf{19.25} & 19.30 & 19.26 & 19.32 \\
				\hdashline
			50    &       & \textbf{22.61} & 22.67 & 22.62 & 22.68 &       & 18.13 & 18.17 & \textbf{17.83} & 17.87 &       & \textbf{23.07} & 23.13 & \textbf{23.07} & 23.13 &       & 17.76 & 17.81 & \textbf{17.74} & 17.78 \\
			&       & \textbf{19.73} & 19.77 & \textbf{19.73} & 19.77 &       & \textbf{21.59} & 21.65 & 21.60 & 21.65 &       & \textbf{23.14} & 23.20 & 23.19 & 23.26 &       & \textbf{18.75} & 18.90 & 18.81 & 18.86 \\
			&       & \textbf{25.06} & 25.11 & 25.07 & 25.13 &       & 24.14 & 24.19 & \textbf{24.12} & 24.18 &       & \textbf{21.65} & 21.70 & 22.03 & 22.08 &       & \textbf{20.59} & 20.64 & 20.78 & 20.84 \\
			&       & \textbf{26.52} & 26.58 & \textbf{26.52} & 26.59 &       & 24.68 & 24.74 & \textbf{24.50} & 24.57 &       & \textbf{17.64} & 17.69 & \textbf{17.64} & 17.69 &       & \textbf{20.47} & 20.53 & 20.91 & 20.96 \\
			&       & \textbf{21.53} & 21.59 & \textbf{21.53} & 21.59 &       & \textbf{24.25} & 24.30 & 24.91 & 24.98 &       & \textbf{20.59} & 20.64 & 20.61 & 20.67 &       & \textbf{22.44} & 22.49 & \textbf{22.44} & 22.51 \\
				\hdashline
			100   &       & \textbf{20.64} & 20.70 & 20.70 & 20.76 &       & \textbf{26.82} & 26.90 & 27.05 & 27.13 &       & 24.67 & 24.73 & \textbf{24.50} & 24.56 &       & 18.66 & 18.70 & \textbf{18.65} & 18.69 \\
			&       & \textbf{21.18} & 21.24 & \textbf{21.18} & 21.24 &       & \textbf{22.01} & 22.07 & \textbf{22.01} & 22.10 &       & \textbf{19.92} & 19.97 & \textbf{19.92} & 19.98 &       & \textbf{24.61} & 24.68 & \textbf{24.61} & 24.68 \\
			&       & \textbf{26.45} & 26.52 & \textbf{26.45} & 26.52 &       & \textbf{21.20} & 21.26 & 21.23 & 21.29 &       & \textbf{26.12} & 26.18 & \textbf{26.12} & 26.19 &       & \textbf{19.68} & 19.73 & 19.89 & 19.95 \\
			&       & \textbf{21.69} & 21.75 & \textbf{21.69} & 21.75 &       & \textbf{25.41} & 25.48 & 25.42 & 25.49 &       & \textbf{21.20} & 21.25 & 21.21 & 21.27 &       & \textbf{22.76} & 22.82 & 22.95 & 23.01 \\
			&       & \textbf{25.70} & 25.76 & \textbf{25.70} & 25.76 &       & 22.65 & 22.70 & \textbf{22.61} & 22.67 &       & \textbf{24.76} & 24.83 & 24.82 & 24.89 &       & \textbf{21.36} & 21.42 & 21.37 & 21.43 \\
			\bottomrule
		\end{tabular}%
	\end{adjustbox}
	\label{tablbub1}%
\end{sidewaystable}%

\begin{sidewaystable}[htbp]
\centering
\caption{Lagrangian relaxation based procedure: percentage gaps of the upper bound w.r.t. the optimal value of the problem for $T=3/S=5$, $T=S=5$, $T=5/S=10$ and $T=S=10$.}
\renewcommand{\arraystretch}{1}
\begin{adjustbox}{max width=1.10\textwidth}
	\begin{tabular}{rrrrrrrrrrrrrrrrrrrrr}
				\toprule
	\multicolumn{21}{c}{Gap$_{UB/{  OPT}}$} \\
	&       & \multicolumn{4}{c}{T=3/S=5}   &       & \multicolumn{4}{c}{T=5/S=5}   &       & \multicolumn{4}{c}{T=5/S=10}  &       & \multicolumn{4}{c}{T=10/S=10} \\
	m/n   &       & Variant 1.iii & Variant 1.iv & Variant 2.iii & Variant 2.iv &       & Variant 1.iii & Variant 1.iv & Variant 2.iii & Variant 2.iv &       & Variant 1.iii & Variant 1.iv & Variant 2.iii & Variant 2.iv &       & Variant 1.iii & Variant 1.iv & Variant 2.iii & Variant 2.iv \\
	\cmidrule{1-1}	\cmidrule{3-6}	\cmidrule{8-11}	\cmidrule{13-16} 	\cmidrule{18-21}
		5     &       & \textbf{$<$0.01} & \textbf{$<$0.01} & \textbf{$<$0.01} & \textbf{$<$0.01} &       & \textbf{$<$0.01} & \textbf{$<$0.01} & \textbf{$<$0.01} & \textbf{$<$0.01} &       & \textbf{$<$0.01} & \textbf{$<$0.01} & \textbf{$<$0.01} & \textbf{$<$0.01} &       & \textbf{$<$0.01} & \textbf{$<$0.01} & \textbf{$<$0.01} & \textbf{$<$0.01} \\
		&       & \textbf{$<$0.01} & \textbf{$<$0.01} & \textbf{$<$0.01} & \textbf{$<$0.01} &       & \textbf{$<$0.01} & \textbf{$<$0.01} & \textbf{$<$0.01} & \textbf{$<$0.01} &       & \textbf{$<$0.01} & \textbf{$<$0.01} & \textbf{$<$0.01} & \textbf{$<$0.01} &       & \textbf{$<$0.01} & \textbf{$<$0.01} & \textbf{$<$0.01} & \textbf{$<$0.01} \\
		&       & \textbf{$<$0.01} & \textbf{$<$0.01} & \textbf{$<$0.01} & \textbf{$<$0.01} &       & \textbf{$<$0.01} & \textbf{$<$0.01} & \textbf{$<$0.01} & \textbf{$<$0.01} &       & \textbf{$<$0.01} & \textbf{$<$0.01} & \textbf{$<$0.01} & \textbf{$<$0.01} &       & \textbf{$<$0.01} & \textbf{$<$0.01} & \textbf{$<$0.01} & \textbf{$<$0.01} \\
		&       & \textbf{$<$0.01} & \textbf{$<$0.01} & \textbf{$<$0.01} & \textbf{$<$0.01} &       & \textbf{$<$0.01} & \textbf{$<$0.01} & \textbf{$<$0.01} & \textbf{$<$0.01} &       & \textbf{$<$0.01} & \textbf{$<$0.01} & \textbf{$<$0.01} & \textbf{$<$0.01} &       & \textbf{$<$0.01} & \textbf{$<$0.01} & \textbf{$<$0.01} & \textbf{$<$0.01} \\
		&       & \textbf{$<$0.01} & \textbf{$<$0.01} & \textbf{$<$0.01} & \textbf{$<$0.01} &       & \textbf{$<$0.01} & \textbf{$<$0.01} & \textbf{$<$0.01} & \textbf{$<$0.01} &       & \textbf{$<$0.01} & \textbf{$<$0.01} & \textbf{$<$0.01} & \textbf{$<$0.01} &       & \textbf{$<$0.01} & \textbf{$<$0.01} & \textbf{$<$0.01} & \textbf{$<$0.01} \\
		\hdashline
		10    &       & \textbf{$<$0.01} & \textbf{$<$0.01} & \textbf{$<$0.01} & \textbf{$<$0.01} &       & \textbf{$<$0.01} & \textbf{$<$0.01} & \textbf{$<$0.01} & \textbf{$<$0.01} &       & \textbf{$<$0.01} & \textbf{$<$0.01} & \textbf{$<$0.01} & \textbf{$<$0.01} &       & 4.86  & 4.86  & \textbf{2.65} & \textbf{2.65} \\
		&       & \textbf{0.94} & \textbf{0.94} & \textbf{0.94} & \textbf{0.94} &       & \textbf{1.93} & \textbf{1.93} & \textbf{1.93} & \textbf{1.93} &       & \textbf{$<$0.01} & \textbf{$<$0.01} & \textbf{$<$0.01} & \textbf{$<$0.01} &       & \textbf{$<$0.01} & \textbf{$<$0.01} & \textbf{$<$0.01} & \textbf{$<$0.01} \\
		&       & \textbf{0.23} & \textbf{0.23} & \textbf{0.23} & \textbf{0.23} &       & \textbf{$<$0.01} & \textbf{$<$0.01} & \textbf{$<$0.01} & \textbf{$<$0.01} &       & \textbf{$<$0.01} & \textbf{$<$0.01} & \textbf{$<$0.01} & \textbf{$<$0.01} &       & \textbf{$<$0.01} & \textbf{$<$0.01} & \textbf{$<$0.01} & \textbf{$<$0.01} \\
		&       & \textbf{$<$0.01} & \textbf{$<$0.01} & \textbf{$<$0.01} & \textbf{$<$0.01} &       & \textbf{$<$0.01} & \textbf{$<$0.01} & \textbf{$<$0.01} & \textbf{$<$0.01} &       & \textbf{$<$0.01} & \textbf{$<$0.01} & \textbf{$<$0.01} & \textbf{$<$0.01} &       & \textbf{$<$0.01} & \textbf{$<$0.01} & \textbf{$<$0.01} & \textbf{$<$0.01} \\
		&       & \textbf{$<$0.01} & \textbf{$<$0.01} & \textbf{$<$0.01} & \textbf{$<$0.01} &       & \textbf{0.17} & \textbf{0.17} & \textbf{0.17} & \textbf{0.17} &       & \textbf{1.99} & \textbf{1.99} & \textbf{1.99} & \textbf{1.99} &       & \textbf{2.44} & \textbf{2.44} & \textbf{2.44} & \textbf{2.44} \\
		\hdashline
		30    &       & \textbf{1.80} & \textbf{1.80} & \textbf{1.80} & \textbf{1.80} &       & \textbf{1.79} & \textbf{1.79} & 2.20  & 2.20  &       & \textbf{$<$0.01} & \textbf{$<$0.01} & \textbf{$<$0.01} & \textbf{$<$0.01} &       & \textbf{0.61} & \textbf{0.61} & 0.67  & 0.67 \\
		&       & \textbf{$<$0.01} & \textbf{$<$0.01} & \textbf{$<$0.01} & \textbf{$<$0.01} &       & \textbf{2.02} & \textbf{2.02} & \textbf{2.02} & \textbf{2.02} &       & \textbf{1.88} & \textbf{1.88} & \textbf{1.88} & \textbf{1.88} &       & \textbf{0.19} & \textbf{0.19} & \textbf{0.19} & 0.57 \\
		&       & \textbf{$<$0.01} & \textbf{$<$0.01} & \textbf{$<$0.01} & \textbf{$<$0.01} &       & \textbf{1.17} & \textbf{1.17} & \textbf{1.17} & \textbf{1.17} &       & \textbf{$<$0.01} & \textbf{$<$0.01} & \textbf{$<$0.01} & \textbf{$<$0.01} &       & \textbf{0.99} & \textbf{0.99} & \textbf{0.99} & \textbf{0.99} \\
		&       & \textbf{0.46} & \textbf{0.46} & \textbf{0.46} & \textbf{0.46} &       & \textbf{1.47} & \textbf{1.47} & \textbf{1.47} & \textbf{1.47} &       & \textbf{0.12} & \textbf{0.12} & \textbf{0.12} & \textbf{0.12} &       & \textbf{0.19} & \textbf{0.19} & \textbf{0.19} & \textbf{0.19} \\
		&       & \textbf{0.41} & \textbf{0.41} & \textbf{0.41} & \textbf{0.41} &       & \textbf{0.75} & \textbf{0.75} & \textbf{0.75} & \textbf{0.75} &       & \textbf{4.59} & \textbf{4.59} & \textbf{4.59} & \textbf{4.59} &       & \textbf{-0.41} & \textbf{-0.41} & \textbf{-0.41} & \textbf{-0.41} \\
		\hdashline
		50    &       & \textbf{0.47} & \textbf{0.47} & \textbf{0.47} & \textbf{0.47} &       & 0.99  & 0.99  & \textbf{0.61} & \textbf{0.61} &       & \textbf{-1.75} & \textbf{-1.75} & \textbf{-1.75} & \textbf{-1.75} &       & -4.70 & -4.70 & \textbf{-4.76} & \textbf{-4.76} \\
		&       & \textbf{$<$0.01} & \textbf{$<$0.01} & \textbf{$<$0.01} & \textbf{$<$0.01} &       & \textbf{0.81} & \textbf{0.81} & \textbf{0.81} & \textbf{0.81} &       & \textbf{-1.13} & \textbf{-1.13} & -1.08 & -1.08 &       & -0.31 & -0.18 & \textbf{-0.31} & \textbf{-0.31} \\
		&       & \textbf{1.47} & \textbf{1.47} & \textbf{1.47} & \textbf{1.47} &       & 1.93  & 1.93  & \textbf{1.89} & \textbf{1.89} &       & \textbf{-3.40} & \textbf{-3.40} & -2.90 & -2.90 &       & \textbf{-4.76} & \textbf{-4.76} & -4.53 & -4.53 \\
		&       & \textbf{3.21} & \textbf{3.21} & \textbf{3.21} & \textbf{3.21} &       & 0.01  & 0.01  & \textbf{-0.23} & \textbf{-0.23} &       & \textbf{0.67} & \textbf{0.67} & \textbf{0.67} & \textbf{0.67} &       & \textbf{0.31} & \textbf{0.31} & 0.78  & 0.78 \\
		&       & \textbf{0.80} & \textbf{0.80} & \textbf{0.80} & \textbf{0.80} &       & \textbf{0.46} & \textbf{0.46} & 1.28  & 1.28  &       & \textbf{1.20} & \textbf{1.20} & \textbf{1.20} & \textbf{1.20} &       & \textbf{-8.42} & \textbf{-8.42} & \textbf{-8.42} & \textbf{-8.42} \\
		\hdashline
		100   &       & \textbf{0.95} & \textbf{0.95} & 0.99  & 0.99  &       & \textbf{-10.39} & \textbf{-10.39} & -10.06 & -10.06 &       & -16.73 & -16.73 & \textbf{-17.00} & \textbf{-17.00} &       & -3.16 & -3.16 & \textbf{-3.19} & \textbf{-3.19} \\
		&       & \textbf{-4.62} & \textbf{-4.62} & \textbf{-4.62} & \textbf{-4.62} &       & \textbf{-19.43} & \textbf{-19.43} & \textbf{-19.43} & -19.39 &       & \textbf{-12.04} & \textbf{-12.04} & \textbf{-12.04} & \textbf{-12.04} &       & \textbf{-18.47} & \textbf{-18.47} & \textbf{-18.47} & \textbf{-18.47} \\
		&       & \textbf{-4.08} & \textbf{-4.08} & \textbf{-4.08} & \textbf{-4.08} &       & \textbf{10.75} & \textbf{10.75} & \textbf{10.75} & \textbf{10.75} &       & \textbf{-4.32} & \textbf{-4.32} & \textbf{-4.32} & \textbf{-4.32} &       & \textbf{-5.56} & \textbf{-5.56} & -5.30 & -5.30 \\
		&       & \textbf{-0.04} & \textbf{-0.04} & \textbf{-0.04} & \textbf{-0.04} &       & \textbf{-14.16} & \textbf{-14.16} & \textbf{-14.16} & \textbf{-14.16} &       & \textbf{-9.70} & \textbf{-9.70} & \textbf{-9.70} & \textbf{-9.70} &       & \textbf{-11.15} & \textbf{-11.15} & -10.89 & -10.89 \\
		&       & \textbf{-1.63} & -1.63 & \textbf{-1.70} & \textbf{-1.70} &       & -8.44 & -8.44 & \textbf{-8.50} & \textbf{-8.50} &       & \textbf{-13.94} & \textbf{-13.94} & -13.87 & -13.87 &       & \textbf{-18.04} & \textbf{-18.04} & \textbf{-18.04} & \textbf{-18.04} \\
		\bottomrule
	\end{tabular}%
\end{adjustbox}
\label{tabubopt1}%
\end{sidewaystable}%

\begin{sidewaystable}[htbp]
\centering
\caption{Lagrangian relaxation versus linear programming relaxation: comparison of lower bounds for $T=3/S=5$, $T=S=5$, $T=5/S=10$ and $T=S=10$.}
\renewcommand{\arraystretch}{1}
\begin{adjustbox}{max width=1.10\textwidth}
	\begin{tabular}{rrrrrrrrrrrrrrrrrrrrr}
					\toprule
	\multicolumn{21}{c}{Gap$_{LP/{  LB}}$} \\
	&       & \multicolumn{4}{c}{T=3/S=5}   &       & \multicolumn{4}{c}{T=5/S=5}   &       & \multicolumn{4}{c}{T=5/S=10}  &       & \multicolumn{4}{c}{T=10/S=10} \\
	m/n   &       & Variant 1.iii & Variant 1.iv & Variant 2.iii & Variant 2.iv &       & Variant 1.iii & Variant 1.iv & Variant 2.iii & Variant 2.iv &       & Variant 1.iii & Variant 1.iv & Variant 2.iii & Variant 2.iv &       & Variant 1.iii & Variant 1.iv & Variant 2.iii & Variant 2.iv \\
	\cmidrule{1-1}	\cmidrule{3-6}	\cmidrule{8-11}	\cmidrule{13-16} 	\cmidrule{18-21}
		5     &       & \textbf{$<$0.01} & 0.03  & \textbf{$<$0.01} & 0.03  &       & \textbf{$<$0.01} & 0.03  & \textbf{$<$0.01} & 0.05  &       & \textbf{$<$0.01} & 0.01  & \textbf{$<$0.01} & 0.01  &       & \textbf{$<$0.01} & 0.08  & \textbf{$<$0.01} & 0.12 \\
		&       & \textbf{$<$0.01} & 0.01  & \textbf{$<$0.01} & 0.01  &       & \textbf{$<$0.01} & 0.04  & \textbf{$<$0.01} & 0.03  &       & \textbf{$<$0.01} & 0.03  & \textbf{$<$0.01} & 0.03  &       & \textbf{$<$0.01} & 0.02  & \textbf{$<$0.01} & 0.03 \\
		&       & \textbf{$<$0.01} & 0.02  & \textbf{$<$0.01} & 0.02  &       & \textbf{$<$0.01} & 0.02  & \textbf{$<$0.01} & 0.02  &       & \textbf{$<$0.01} & $<$0.01  & \textbf{$<$0.01} & $<$0.01  &       & \textbf{$<$0.01} & 0.02  & \textbf{$<$0.01} & 0.02 \\
		&       & \textbf{$<$0.01} & \textbf{$<$0.01} & \textbf{$<$0.01} & $<$0.01  &       & \textbf{0.01} & \textbf{0.01} & \textbf{0.01} & \textbf{0.01} &       & \textbf{$<$0.01} & 0.04  & \textbf{$<$0.01} & 0.03  &       & \textbf{$<$0.01} & 0.04  & \textbf{$<$0.01} & 0.04 \\
		&       & \textbf{$<$0.01} & 0.02  & \textbf{$<$0.01} & 0.02  &       & \textbf{0.01} & \textbf{0.01} & \textbf{0.01} & \textbf{0.01} &       & \textbf{$<$0.01} & 0.01  & \textbf{$<$0.01} & 0.03  &       & \textbf{$<$0.01} & 0.01  & \textbf{$<$0.01} & 0.03 \\
		\hdashline
		10    &       & \textbf{$<$0.01} & 0.05  & \textbf{$<$0.01} & 0.06  &       & \textbf{$<$0.01} & 0.01  & \textbf{$<$0.01} & 0.01  &       & \textbf{$<$0.01} & 0.04  & \textbf{$<$0.01} & 0.04  &       & \textbf{$<$0.01} & 0.05  & \textbf{$<$0.01} & 0.04 \\
		&       & \textbf{$<$0.01} & 0.02  & \textbf{$<$0.01} & 0.04  &       & \textbf{$<$0.01} & 0.04  & \textbf{$<$0.01} & 0.06  &       & \textbf{$<$0.01} & 0.03  & \textbf{$<$0.01} & 0.04  &       & \textbf{$<$0.01} & 0.02  & \textbf{$<$0.01} & 0.03 \\
		&       & \textbf{$<$0.01} & 0.03  & \textbf{$<$0.01} & 0.02  &       & \textbf{$<$0.01} & 0.03  & \textbf{$<$0.01} & 0.03  &       & \textbf{$<$0.01} & 0.03  & \textbf{$<$0.01} & 0.05  &       & \textbf{$<$0.01} & 0.04  & \textbf{$<$0.01} & 0.06 \\
		&       & \textbf{$<$0.01} & 0.01  & \textbf{$<$0.01} & 0.02  &       & \textbf{$<$0.01} & 0.02  & \textbf{$<$0.01} & 0.03  &       & \textbf{$<$0.01} & 0.03  & \textbf{$<$0.01} & 0.02  &       & \textbf{$<$0.01} & 0.07  & \textbf{$<$0.01} & 0.05 \\
		&       & \textbf{$<$0.01} & 0.06  & \textbf{$<$0.01} & 0.08  &       & \textbf{$<$0.01} & 0.04  & \textbf{$<$0.01} & 0.05  &       & \textbf{$<$0.01} & 0.07  & \textbf{$<$0.01} & 0.09  &       & \textbf{$<$0.01} & 0.05  & \textbf{$<$0.01} & 0.08 \\
		\hdashline
		30    &       & \textbf{$<$0.01} & 0.08  & \textbf{$<$0.01} & 0.09  &       & \textbf{$<$0.01} & 0.06  & 0.16  & 0.23  &       & \textbf{$<$0.01} & 0.05  & \textbf{$<$0.01} & 0.05  &       & \textbf{0.02} & 0.07  & \textbf{0.04} & 0.09 \\
		&       & \textbf{$<$0.01} & 0.06  & \textbf{$<$0.01} & 0.07  &       & \textbf{$<$0.01} & 0.06  & $<$0.01  & 0.08  &       & \textbf{$<$0.01} & 0.08  & \textbf{$<$0.01} & 0.08  &       & \textbf{$<$0.01} & 0.06  & \textbf{$<$0.01} & 0.06 \\
		&       & \textbf{$<$0.01} & 0.06  & \textbf{$<$0.01} & 0.07  &       & \textbf{$<$0.01} & 0.07  & 0.02  & 0.09  &       & \textbf{$<$0.01} & 0.06  & \textbf{$<$0.01} & 0.07  &       & \textbf{$<$0.01} & 0.06  & \textbf{$<$0.01} & 0.06 \\
		&       & \textbf{$<$0.01} & 0.08  & \textbf{$<$0.01} & 0.09  &       & \textbf{$<$0.01} & 0.04  & 0.02  & 0.08  &       & \textbf{$<$0.01} & 0.05  & \textbf{$<$0.01} & 0.06  &       & \textbf{0.02} & 0.07  & 0.05  & 0.09 \\
		&       & \textbf{$<$0.01} & 0.07  & \textbf{$<$0.01} & 0.08  &       & \textbf{$<$0.01} & 0.03  & $<$0.01  & 0.04  &       & \textbf{$<$0.01} & 0.06  & \textbf{$<$0.01} & 0.06  &       & \textbf{$<$0.01} & 0.06  & 0.02  & 0.08 \\
		\hdashline
		50    &       & \textbf{$<$0.01} & 0.08  & 0.02  & 0.09  &       & \textbf{$<$0.01} & 0.06  & 0.02  & 0.07  &       & \textbf{$<$0.01} & 0.08  & \textbf{$<$0.01} & 0.08  &       & \textbf{0.01} & 0.06  & 0.02  & 0.08 \\
		&       & \textbf{$<$0.01} & 0.06  & \textbf{$<$0.01} & 0.06  &       & \textbf{$<$0.01} & 0.07  & 0.02  & 0.08  &       & 0.01  & 0.09  & 0.02  & 0.11  &       & \textbf{0.08} & 0.13  & 0.15  & 0.21 \\
		&       & \textbf{0.01} & 0.09  & 0.02  & 0.11  &       & \textbf{$<$0.01} & 0.07  & \textbf{$<$0.01} & 0.08  &       & \textbf{$<$0.01} & 0.07  & \textbf{$<$0.01} & 0.07  &       & \textbf{0.01} & 0.07  & 0.04  & 0.11 \\
		&       & \textbf{$<$0.01} & 0.09  & \textbf{$<$0.01} & 0.10  &       & \textbf{$<$0.01} & 0.08  & \textbf{$<$0.01} & 0.09  &       & \textbf{$<$0.01} & 0.05  & \textbf{$<$0.01} & 0.06  &       & \textbf{$<$0.01} & 0.07  & 0.08  & 0.14 \\
		&       & \textbf{$<$0.01} & 0.07  & \textbf{$<$0.01} & 0.07  &       & \textbf{0.01} & 0.07  & 0.07  & 0.15  &       & \textbf{0.02} & 0.08  & 0.04  & 0.11  &       & \textbf{$<$0.01} & 0.07  & 0.01  & 0.09 \\
		\hdashline
		100   &       & \textbf{0.04} & 0.11  & 0.07  & 0.15  &       & \textbf{$<$0.01} & 0.10  & 0.01  & 0.12  &       & \textbf{0.01} & 0.09  & 0.02  & 0.10  &       & \textbf{$<$0.01} & 0.06  & 0.01  & 0.07 \\
		&       & \textbf{$<$0.01} & 0.08  & \textbf{$<$0.01} & 0.08  &       & \textbf{0.01} & 0.08  & \textbf{0.01} & 0.10  &       & \textbf{$<$0.01} & 0.07  & 0.01  & 0.07  &       & \textbf{$<$0.01} & 0.09  & 0.01  & 0.10 \\
		&       & \textbf{$<$0.01} & 0.10  & \textbf{$<$0.01} & 0.10  &       & \textbf{0.01} & 0.08  & 0.04  & 0.12  &       & \textbf{0.01} & 0.10  & \textbf{0.01} & 0.11  &       & \textbf{$<$0.01} & 0.07  & 0.02  & 0.09 \\
		&       & \textbf{$<$0.01} & 0.08  & 0.01  & 0.08  &       & \textbf{0.02} & 0.10  & \textbf{0.02} & 0.12  &       & \textbf{0.01} & 0.08  & 0.03  & 0.10  &       & \textbf{0.02} & 0.09  & \textbf{0.02} & 0.10 \\
		&       & \textbf{$<$0.01} & 0.09  & 0.07  & 0.16  &       & \textbf{0.01} & 0.08  & 0.02  & 0.09  &       & \textbf{$<$0.01} & 0.10  & 0.03  & 0.12  &       & \textbf{$<$0.01} & 0.08  & 0.01  & 0.08 \\
		\bottomrule
	\end{tabular}%
\end{adjustbox}
\label{tablplag1}%
\end{sidewaystable}%

\begin{sidewaystable}[htbp]
	\centering
	\caption{ Lagrangian relaxation based procedure: CPU time (seconds) ($T=3,S=5$ and $T=S=5$)}
	\renewcommand{\arraystretch}{0.95}
	\begin{adjustbox}{max width=1.00\textwidth}
	\begin{tabular}{rrrrrrrrrrrrrrrrr}
		\toprule
		\multicolumn{17}{c}{Lagrangian relaxation times} \\
		m/n   &       & T/S   &       & Variant 1.iii & Variant 1.iv & Variant 2.iii & Variant 2.iv & $t_{solv}$ &       &       &       & Variant 1.iii & Variant 1.iv & Variant 2.iii & Variant 2.iv & $t_{solv}$ \\
		\midrule
		5     &       & 3/5   &       & 1.90  & 0.67  & 2.23  & \textbf{0.88}  & 0.22  &       & 5/5   &       & 2.48  & 0.97  & 2.45  & \textbf{0.51}  & 0.24 \\
		&       &       &       & 1.89  & 1.08  & 1.98  &\textbf{0.42}  & 0.00  &       &       &       & 2.40  & 1.11  & 2.45  & \textbf{0.41}  & 0.01 \\
		&       &       &       & 1.84  & 0.58  & 1.93  & \textbf{0.53}  & 0.02  &       &       &       & 2.42  & 0.72  & 2.39  & \textbf{0.34}  & 0.01 \\
		&       &       &       & 1.81  & 0.38  & 1.78  & \textbf{0.22}  & 0.00  &       &       &       & 0.08  & 0.11  & 0.09  & \textbf{0.09}  & 0.00 \\
		&       &       &       & 1.90  & 0.62  & 1.79  &\textbf{ 0.31}  & 0.01  &       &       &       & 0.16  & 0.16  & 0.16  & \textbf{0.16}  & 0.00 \\
		\hdashline
		10    &       & 3/5   &       & 3.07  & 1.08  & 3.04  &\textbf{ 0.70}  & 0.04  &       & 5/5   &       & 6.27  & 2.51  & 6.04  & \textbf{1.28}  & 0.03 \\
		&       &       &       & 3.14  & 1.26  & 3.03  &\textbf{ 0.66}  & 0.03  &       &       &       & 6.13  & 2.09  & 6.12  & \textbf{1.30}  & 0.10 \\
		&       &       &       & 3.14  & 1.09  & 3.37  & \textbf{0.82}  & 0.30  &       &       &       & 6.10  & 2.28  & 5.91  & \textbf{1.36}  & 0.09 \\
		&       &       &       & 3.00  & 1.28  & 3.15  & \textbf{0.77}  & 0.02  &       &       &       & 6.19  & 2.25  & 6.13  & \textbf{1.53}  & 0.05 \\
		&       &       &       & 2.96  & 1.20  & 3.35  & \textbf{0.71}  & 0.03  &       &       &       & 6.32  & 2.28  & 5.98  & \textbf{1.39}  & 0.03 \\
		\hdashline
	30    &       & 3/5   &       & 16.86 & 8.25  & 16.49 &\textbf{ 4.27}  & 2.66  &       & 5/5   &       & 34.21 & 16.58 & 33.77 & \textbf{8.16}  & 30.26 \\
		&       &       &       & 16.29 & 6.37  & 17.02 & \textbf{4.06}  & 14.33 &       &       &       & 33.93 & 14.41 & 33.42 & \textbf{8.24}  & 3274.71 \\
		&       &       &       & 16.04 & 7.46  & 15.73 & \textbf{4.29}  & 1.84  &       &       &       & 34.02 & 14.27 & 33.38 & \textbf{7.83}  & 973.39 \\
		&       &       &       & 16.26 & 5.87  & 16.01 & \textbf{4.17}  & 25.74 &       &       &       & 34.20 & 15.96 & 33.57 & \textbf{8.42}  & 58.44 \\
		&       &       &       & 16.33 & 6.68  & 15.90 & \textbf{4.69}  & 8.38  &       &       &       & 34.20 & 16.21 & 33.51 & 9.20  & 24.62 \\
		\hdashline
		50    &       & 3/5   &       & 38.00 & 15.32 & 37.50 & \textbf{9.69}  & 573.93 &       & 5/5   &       & 87.70 & 42.51 & 86.58 & \textbf{25.63} & 10803.50* \\
		&       &       &       & 41.18 & 19.73 & 40.54 & \textbf{13.27} & 12.07 &       &       &       & 88.17 & 42.23 & 86.92 & \textbf{25.74} & 10805.10* \\
		&       &       &       & 39.84 & 15.79 & 39.13 & \textbf{9.11}  & 1810.31 &       &       &       & 88.48 & 35.55 & 87.27 & \textbf{26.63} & 10804.40* \\
		&       &       &       & 40.28 & 15.07 & 39.42 & \textbf{9.82}  & 9295.68 &       &       &       & 90.25 & 34.55 & 88.97 & \textbf{22.99} & 10806.50* \\
		&       &       &       & 39.22 & 15.40 & 38.66 & \textbf{10.83} & 372.20 &       &       &       & 87.39 & 42.79 & 86.03 & \textbf{23.35} & 10806.60* \\
			\hdashline
				100   &       & 3/5   &       & 142.07 & 52.81 & 140.39 & \textbf{44.95} & 10807.90* &       & 5/5   &       & 350.30 & 192.99 & 346.24 & \textbf{96.28} & 10837.40* \\
		&       &       &       & 148.58 & 57.21 & 146.94 & \textbf{48.72} & 10808.90* &       &       &       & 350.97 & 144.36 & 346.38 & \textbf{104.88} & 10809.50* \\
		&       &       &       & 149.82 & 63.96 & 147.86 & \textbf{53.27} & 10807.90* &       &       &       & 348.24 & 157.14 & 343.69 & \textbf{90.15} & Out of mem.* \\
		&       &       &       & 146.31 & 66.77 & 144.21 & \textbf{45.52} & 10808.50* &       &       &       & 344.73 & 127.22 & 340.27 & \textbf{101.62} & 10811.00* \\
		&       &       &       & 149.01 & 74.91 & 147.30 & \textbf{38.27} & 10808.70* &       &       &       & 343.20 & 164.47 & 339.04 & \textbf{112.76} & 10809.40* \\
		\bottomrule
	\end{tabular}%
\end{adjustbox}
\label{tabtimes1}%
\end{sidewaystable}%

\begin{sidewaystable}[htbp]
	\centering
\caption{ Lagrangian relaxation based procedure: CPU time (seconds) ($T=5,S=10$ and $T=S=10$)}
\renewcommand{\arraystretch}{0.95}
\begin{adjustbox}{max width=1.00\textwidth}
	\begin{tabular}{rrrrrrrrrrrrrrrrr}
		\toprule
		\multicolumn{17}{c}{Lagrangian relaxation times} \\
		m/n   &       & T/S   &       & Variant 1.iii & Variant 1.iv & Variant 2.iii & Variant 2.iv & $t_{solv}$ &       &       &       & Variant 1.iii & Variant 1.iv & Variant 2.iii & Variant 2.iv & $t_{solv}$ \\
		\cmidrule{1-1} \cmidrule{3-3} \cmidrule{5-9} \cmidrule{11-11} \cmidrule{13-17}
		5     &       & 5/10  &       & 3.90  & 0.73  & 4.41  & \textbf{0.58}  & 0.24  &       & 10/10 &       & 8.81  & 2.53  & 8.71  & \textbf{1.59}  & 0.23 \\
		&       &       &       & 4.57  & 0.92  & 3.67  & \textbf{0.64}  & 0.03  &       &       &       & 8.96  & 3.31  & 8.71  & \textbf{1.69}  & 0.08 \\
		&       &       &       & 4.32  & 1.72  & 3.85  & \textbf{0.76}  & 0.01  &       &       &       & 8.77  & 2.31  & 8.49  &\textbf{ 1.42}  & 0.02 \\
		&       &       &       & 4.20  & 1.98  & 3.60  &\textbf{ 0.67}  & 0.03  &       &       &       & 8.89  & 3.37  & 8.63  & \textbf{2.22}  & 0.01 \\
		&       &       &       & 4.23  & 1.40  & 3.84  & \textbf{0.87}  & 0.01  &       &       &       & 9.05  & 2.67  & 8.57  & \textbf{1.61}  & 0.30 \\
			\hdashline
		10    &       & 5/10  &       & 9.64  & 3.09  & 9.48  & \textbf{1.75}  & 0.12  &       & 10/10 &       & 27.39 & 11.72 & 26.82 & \textbf{6.10}  & 1.35 \\
		&       &       &       & 9.95  & 5.26  & 9.78  & \textbf{2.93}  & 0.43  &       &       &       & 27.00 & 9.39  & 26.33 & \textbf{5.90}  & 0.27 \\
		&       &       &       & 9.74  & 4.23  & 9.59  & \textbf{2.18}  & 0.37  &       &       &       & 27.41 & 9.53  & 26.97 & \textbf{5.57}  & 1.07 \\
		&       &       &       & 9.72  & 3.34  & 9.42  & \textbf{1.76}  & 0.04  &       &       &       & 27.13 & 10.28 & 26.57 & \textbf{5.34 } & 0.31 \\
		&       &       &       & 9.78  & 4.03  & 9.53  & \textbf{2.00}  & 0.49  &       &       &       & 26.99 & 7.79  & 26.49 & \textbf{6.15}  & 0.29 \\
		\hdashline
		30    &       & 5/10  &       & 63.40 & 28.17 & 62.67 & \textbf{17.00} & 943.88 &       & 10/10 &       & 216.05 & 104.32 & 213.07 & \textbf{57.24} & 10805.10* \\
		&       &       &       & 64.44 & 31.75 & 63.41 & \textbf{18.85} & 10805.90* &       &       &       & 215.70 & 101.18 & 212.63 & \textbf{57.52} & 10801.90* \\
		&       &       &       & 64.21 & 31.40 & 63.27 & \textbf{16.71} & 1242.70 &       &       &       & 216.90 & 103.88 & 217.32 & \textbf{49.34} & 10805.30* \\
		&       &       &       & 63.76 & 30.94 & 62.92 & \textbf{16.24} & 10805.40* &       &       &       & 214.30 & 97.44 & 211.24 & \textbf{71.92} & 10805.00* \\
		&       &       &       & 64.12 & 25.69 & 63.26 & \textbf{15.87} & 10804.70* &       &       &       & 215.20 & 93.24 & 212.24 & \textbf{46.78} & 10805.30* \\
		\hdashline
	50    &       & 5/10  &       & 175.28 & 75.15 & 172.85 & \textbf{45.71} & 10808.30* &       & 10/10 &       & 592.43 & 241.40 & 584.85 & \textbf{168.82} & 10808.20* \\
		&       &       &       & 170.21 & 83.57 & 167.98 & \textbf{46.79} & 10818.80* &       &       &       & 595.36 & 278.10 & 588.00 & \textbf{171.96} & 10806.90* \\
		&       &       &       & 172.01 & 71.11 & 169.95 & \textbf{54.46} & 10806.80* &       &       &       & 590.74 & 293.53 & 583.05 & \textbf{182.19} & 10807.10* \\
		&       &       &       & 176.84 & 61.79 & 174.55 & \textbf{48.25} & 3188.16 &       &       &       & 593.77 & 237.06 & 586.50 & \textbf{164.41} & 10810.40* \\
		&       &       &       & 172.13 & 89.65 & 169.95 & \textbf{50.64} & 10807.60* &       &       &       & 593.08 & 291.27 & 585.22 & \textbf{153.64} & 10859.20* \\
		\hdashline
		100   &       & 5/10  &       & 679.40 & 331.05 & 670.86 & \textbf{208.37} & 10827.10* &       & 10/10 &       & 2351.86 & 1187.46 & 2320.97 & \textbf{814.57} & 10845.60* \\
		&       &       &       & 681.25 & 311.92 & 672.81 & \textbf{249.65} & 10810.30* &       &       &       & 2362.42 & 1272.60 & 2330.85 & \textbf{812.97} & 10810.00* \\
		&       &       &       & 679.43 & 328.07 & 671.00 & \textbf{189.48} & 10807.30* &       &       &       & 2359.49 & 1163.45 & 2329.60 &\textbf{ 672.58} & 10828.40* \\
		&       &       &       & 681.71 & 336.42 & 673.59 &\textbf{ 219.87} & 10807.50* &       &       &       & 2358.36 & 1467.57 & 2328.16 & \textbf{784.03} & 10810.60* \\
		&       &       &       & 688.06 & 293.90 & 679.29 & \textbf{212.49} & 10808.70* &       &       &       & 2355.81 & 959.28 & 2324.92 & \textbf{675.98} & 10811.80* \\
		\bottomrule
	\end{tabular}%
\end{adjustbox}
\label{tabtimes2}%
\end{sidewaystable}%

{
	Finally, in Table \ref{tab1}, we have compared LB and UB of Lagrangian relaxation with the best solution and lower bound provided by CPLEX using the default parameters (Balanced column) and by emphasizing feasibility (Feasible column). In this table we have reported the average gaps using variant 2.iv of the Lagrangian based heuristic because this variant provides the best running times and it is among the ones with best gaps. Observe that the difference of the total average using balanced or feasible option is not larger than $1.1\%$ (absolute value) in each of the four analyzed gaps, i.e., not big differences are appreciated by these changes of parameters of CPLEX. On the other hand, we remark the high improvement of the lower bound of CPLEX with respect to the value of the Lagrangian relaxation (very close to the linear relaxation, see Table 4).
\begin{table}[htbp]
	\centering
	\caption{{Table including average gaps between $UB$ and $LB$ of variant 2.iv Lagrangian based heuristic and best lower bound ($LB_{cplex}$) and best solution (OPT) using CPLEX.}}
	\scalebox{0.75}{
		{ 
		\begin{tabular}{l|rr|rr|rr|rr}
			& \multicolumn{2}{l|}{{\textbf{Gap $LB$ and $LB_{cplex}$} }}        & \multicolumn{2}{l|}{{\textbf{Gap $UB$ and $LB_{cplex}$} }}      & \multicolumn{2}{l|}{{\textbf{Gap $LB$ and $OPT$ }}}        & \multicolumn{2}{l}{{\textbf{Gap $UB$ and $OPT$ }}}   \\
			& \multicolumn{1}{l}{\textbf{Balanced}} & \multicolumn{1}{l|}{\textbf{Feasible}} & \multicolumn{1}{l}{\textbf{Balanced}} & \multicolumn{1}{l|}{\textbf{Feasible}} & \multicolumn{1}{l}{\textbf{Balanced}} & \multicolumn{1}{l|}{\textbf{Feasible}} & \multicolumn{1}{l}{\textbf{Balanced}} & \multicolumn{1}{l}{\textbf{Feasible}} \\
			\hline
			\textbf{(T,S)=(3/3)}&&&&&&&&\\
			(m,n)=(5/5)   & 17.58\% & 17.58\% & 0.00\% & 0.00\% & 17.58\% & 17.58\% & 0.00\% & 0.00\% \\
			(m,n)=(10/10) & 21.98\% & 21.98\% & 0.39\% & 0.39\% & 21.98\% & 21.98\% & 0.39\% & 0.39\% \\
			(m,n)=(30/30) & 24.12\% & 24.12\% & 1.57\% & 1.57\% & 24.13\% & 24.13\% & 1.56\% & 1.56\% \\
			(m,n)=(50/50) & 20.50\% & 20.03\% & 1.08\% & 1.65\% & 20.51\% & 20.52\% & 1.06\% & 1.06\% \\
			(m,n)=(100/100)& 19.33\% & 17.14\% & 6.74\% & 9.19\% & 24.74\% & 24.02\% & 0.01\% & 0.98\% \\
			\textbf{Average} & \textbf{20.70\%} & \textbf{20.17\%} & \textbf{1.96\%} & \textbf{2.56\%} & \textbf{21.79\%} & \textbf{21.65\%} & \textbf{0.60\%} & \textbf{0.80\%} \\
			\hline
			\textbf{(T,S)=(3/5)}&&&&&&&&\\
			(m,n)=(5/5)  & 7.70\% & 7.70\% & 0.00\% & 0.00\% & 7.70\% & 7.70\% & 0.00\% & 0.00\% \\
			(m,n)=(10/10) & 19.80\% & 19.80\% & 0.23\% & 0.23\% & 19.80\% & 19.80\% & 0.23\% & 0.23\% \\
			(m,n)=(30/30) & 22.95\% & 22.95\% & 0.55\% & 0.55\% & 22.97\% & 22.97\% & 0.53\% & 0.53\% \\
			(m,n)=(50/50)& 22.22\% & 20.43\% & 1.22\% & 3.37\% & 22.24\% & 22.27\% & 1.19\% & 1.16\% \\
			(m,n)=(100/100) & 16.28\% & 14.70\% & 8.26\% & 9.97\% & 24.57\% & 23.71\% & -1.89\% & -0.68\% \\
			\textbf{Average} & \textbf{17.79\%} & \textbf{17.12\%} & \textbf{2.05\%} & \textbf{2.82\%} & \textbf{19.46\%} & \textbf{19.29\%} & \textbf{0.01\%} & \textbf{0.25\%} \\
			\hline
			\textbf{(T,S)=(5/5)}&&&&&&&&\\
			(m,n)=(5/5)   & 11.21\% & 11.21\% & 0.00\% & 0.00\% & 11.21\% & 11.21\% & 0.00\% & 0.00\% \\
			(m,n)=(10/10) & 14.73\% & 14.73\% & 0.42\% & 0.42\% & 14.73\% & 14.73\% & 0.42\% & 0.42\% \\
			(m,n)=(30/30)& 17.58\% & 17.58\% & 1.54\% & 1.54\% & 17.59\% & 17.59\% & 1.52\% & 1.52\% \\
			(m,n)=(50/50) & 17.18\% & 12.63\% & 6.58\% & 11.42\% & 21.97\% & 22.36\% & 0.87\% & 0.35\% \\
			(m,n)=(100/100) & 14.24\% & 12.41\% & 11.04\% & 12.94\% & 29.31\% & 25.21\% & -8.12\% & -2.01\% \\
			\textbf{Average} & \textbf{14.99\%} & \textbf{13.71\%} & \textbf{3.92\%} & \textbf{5.26\%} & \textbf{18.96\%} & \textbf{18.22\%} & \textbf{-1.06\%} & \textbf{0.06\%} \\
			\hline
			\textbf{(T,S)=(5/10)}&&&&&&&&\\
			(m,n)=(5/5)  & 11.28\% & 11.28\% & 0.00\% & 0.00\% & 11.28\% & 11.28\% & 0.00\% & 0.00\% \\
			(m,n)=(10/10) & 19.06\% & 19.06\% & 0.40\% & 0.40\% & 19.06\% & 19.06\% & 0.40\% & 0.40\% \\
			(m,n)=(30/30)& 18.21\% & 14.52\% & 2.46\% & 6.42\% & 19.16\% & 19.37\% & 1.32\% & 1.05\% \\
			(m,n)=(50/50) & 14.19\% & 11.57\% & 8.29\% & 11.04\% & 21.93\% & 20.96\% & -0.77\% & 0.49\% \\
			(m,n)=(100/100) & 12.27\% & 11.79\% & 12.68\% & 13.15\% & 31.12\% & 26.88\% & -11.38\% & -4.87\% \\
			\textbf{Average} & \textbf{15.00\%} & \textbf{13.64\%} & \textbf{4.77\%} & \textbf{6.20\%} & \textbf{20.51\%} & \textbf{19.51\%} & \textbf{-2.09\%} & \textbf{-0.59\%} \\
			\hline
			\textbf{(T,S)=(10/10)}&&&&&&&&\\
			(m,n)=(5/5)   & 17.25\% & 17.25\% & 0.00\% & 0.00\% & 17.25\% & 17.25\% & 0.00\% & 0.00\% \\
			(m,n)=(10/10) & 20.02\% & 20.02\% & 1.02\% & 1.02\% & 20.02\% & 20.02\% & 1.02\% & 1.02\% \\
			(m,n)=(30/30) & 15.10\% & 10.03\% & 4.11\% & 9.51\% & 18.33\% & 18.57\% & 0.40\% & 0.11\% \\
			(m,n)=(50/50) & 10.65\% & 9.33\% & 10.67\% & 11.97\% & 22.75\% & 21.33\% & -3.45\% & -1.47\% \\
			(m,n)=(100/100) & 9.84\% & 9.71\% & 12.99\% & 13.11\% & 29.12\% & 25.67\% & -11.18\% & -5.58\% \\
			\textbf{Average} & \textbf{14.57\%} & \textbf{13.27\%} & \textbf{5.76\%} & \textbf{7.12\%} & \textbf{21.49\%} & \textbf{20.57\%} & \textbf{-2.64\%} & \textbf{-1.19\%} \\
			\hline
			\hline
			\textbf{Total Av.} & \textbf{16.61\%} & \textbf{15.58\%} & \textbf{3.69\%} & \textbf{4.79\%} & \textbf{20.44\%} & \textbf{19.85\%} & \textbf{-1.03\%} & \textbf{-0.13\%} \\
			\hline
		\end{tabular}%
	}
	}
	\label{tab1}%
\end{table}%

}

\subsection{The relevance of considering a multi-period stochastic setting}\label{EVPI}
\label{subsec:EVPI-VMS}
In this work, we have proposed a general formulation for covering location problems that captures the multi-period and stochastic nature of some problems. 
This led us to a very general and thus potentially more involved model.
One important question that we have formulated and that it is still to be answered is: does it compensate using such a more involved model instead of a simplified one (e.g. static and/or deterministic)?

As far as the relevance of considering a stochastic programming modeling framework for a problem is concerned, we know that there are no ``robust'' measures for evaluating that relevance (see, for instance, \citealt{Birge&Louveaux:2011}). 
However, as explained in \cite{Correia2015}, two measures are often used as an indication of the convenience of an stochastic model. 
These are the value of the stochastic solution (VSS) and the expected value of perfect information (EVPI). 
The former relies on the resolution of the so-called expected value problem which consists of replacing the random variables by their expectation.
The resulting solution, if feasible to the stochastic problem, is then evaluated as an approximate stochastic (first-stage) solution and its value compared with the optimal value to the stochastic problem.
In our case, this cannot be accomplished since we have binary stochastic parameters 
($a_{ijt}^s$) and thus the expected value problem is meaningless. 
Therefore, we resort to the second measure above mentioned, which is the EVPI.
This is a value that indicates the amount that the decision maker should be willing to pay to get access to perfect information about the future.

The EVPI is computed as the difference between the objective value of the stochastic problem (SP) and the wait-and-see solution
(WS). 
The WS is the weighted sum for the optimal values to all single-scenario problems considering as weights the corresponding probabilities for the scenarios. 
Accordingly, $\mbox{WS}=\sum_{s\in {\cal S}}\pi_s \mbox{DP}_s$, where $\mbox{DP}_s$ is the optimal value of the 
stochastic problem related with scenario $s\in {\cal S}$; EVPI $=$ SP $-$ WS.

With respect to the relevance of using a multi-period modeling framework, we consider a measure described by {\cite{Nickel2015}} as the value of the multi-period solution (VMS). 
It consists of comparing the optimal value of the multi-period model (MPS) with the value of a feasible solution (1PS) to that model obtained by solving to optimality a static counterpart.
This way, we are evaluating how good a `static' solution 
may be {when considered} 
as an approximate solution to the multi-period problem. 

In our case, we obtain a static counterpart to our problem by assuming that opening and closing decisions can only be made in the first period. 
Consequently, the static counterpart is obtained from the original model setting 

$z_{it}=0$, for $i\in I$, $t\in {\cal T}\setminus\{1\}$ 
and 
$z^\prime_{it}=0$ for $i\in I$ and $t\in {\cal T}\setminus\{1,T\}$. 

Finally, we have VMS $=$1PS $-$ MPS.

In Table~\ref{tab10} we can observe the results for the above-mentioned measures: EVPI  and VMS.
Then, we can see that in average, the EVPI is $63.98$ and VMS is $323.43$. 
As a conclusion, both measures show that it is important to take into account stochastic and multi-period aspects.

\begin{table}[htbp]
	\centering
\caption{ The relevance of considering a multi-period stochastic setting: EVPI and VMS.}
\renewcommand{\arraystretch}{0.65}
\begin{adjustbox}{max width=1.00\textwidth}
	\begin{tabular}{rrrrrrrrr}
		&&&&&&\\
\toprule
m/n   && $T/S$   && EVPI && VMS &&  SP/MPS\\
\cmidrule{1-1}\cmidrule{3-3}\cmidrule{5-5}\cmidrule{7-7} \cmidrule{9-9}
		5     &       & 3     &       & 4.11  && 0.00 && 15.03 \\
		&       &       &       & 0.00  && 0.00 && -9.24\\
		&       &       &       & 0.29  && 0.00 &&31.96\\
		&       &       &       & 7.08  && 3.25 &&12.79\\
		&       &       &       & 2.33  && 2.24 &&14.19\\
		\hdashline
		10    &       & 3     &       & 14.58 && 37.32 &&31.47\\
		&       &       &       & 5.47  && 51.20 &&-22.34\\
		&       &       &       & 11.99 && 0.00 &&-38.89\\
		&       &       &       & 22.54 && 7.76 &&-13.12\\
		&       &       &       & 2.31  && 6.46 &&-0.04\\
		\hdashline
		30    &       & 3     &       & 42.83 && 616.84 &&-1341.48\\
		&       &       &       & 52.99 && 285.87&& -708.26\\
		&       &       &       & 136.94 && 290.74 &&-889.88\\
		&       &       &       & 69.82 && 254.51&& -1246.03\\
		&       &       &       & 187.12 && 226.15 &&-680.27\\
		\hdashline
		50    &       & 3     &       & 127.40 && 385.67&&-2329.67 \\
		&       &       &       & 147.46 && 1286.12&& -3378.28\\
		&       &       &       & 66.15 && 860.84&& -2914.91\\
		&       &       &       & 204.61 && 1470.84 &&-2404.21\\
		&       &       &       & 173.63 && 683.48 &&-2268.85\\
		\hdashline
		\multicolumn{4}{c}{Average}&63.98&&323.43&&-907.00\\
		\bottomrule
		\end{tabular}
\end{adjustbox}
\label{tab10}%
\end{table}%

\section{Conclusions}
\label{conclusions}

A very general set covering location problem with stochasticity and multiple time periods has been
introduced and analyzed in this paper. 
Almost all the covering location problems in the literature can be 
looked at as particular cases of this general problem, which in consequence is worth studying
but difficult to solve. 
An integer linear programming formulation has been developed which on the one hand is able to 
deal with instances of limited size and, on the other hand, has been the basis for a Lagrangian heuristic designed to obtain good solutions for larger instances. 
Despite the fact that the relaxed Lagrangian subproblem presents the integrality property, the goodness of the solutions it generates in a subgradient {optimization} framework has been deduced from an extensive computational study.

As a matter of future research, additional properties and improvements of the formulation (by means of ad-hoc valid inequalities) and/or some of their possible particularizations will possibly lead to reduced computational times and to the possibility of optimally solving larger instances.

\section*{Acknowledgements}
The research of Alfredo Mar\'{\i}n  is partially funded by  Spanish
{\sl Ministerio de Econom\'{\i}a y Competitividad}, project MTM2015-65915-R,
{\sl  Fundaci\'on S\'eneca de la Consejer\'{\i}a de Educaci\'on de la Comunidad Aut\'onoma de
	la Regi\'on de Murcia}, project 19320/PI/14, and {\sl Fundaci\'on BBVA}, project
``Cost-sensitive classification. A mathematical optimization approach'' (COSECLA).
Luisa I. Mart{\'\i}nez Merino and A.M. Rodr{\'\i}guez-Ch{\'\i}a  acknowledge that research reported here was supported by MTM2013-46962-C2-2-P, MTM2016-74983-C2-2-R, supported by Agencia Estatal de
Investigaci\'on (AEI) and the European Regional Development's funds (FEDER).
Francisco Saldanha-da-Gama was supported by the Portuguese Science Foundation (FCT---Fundação para a Ciência e Tecnologia) under the project UID/MAT/ 04561/2013 (CMAF-CIO/FCUL).
This support is gratefully acknowledged.

\bibliographystyle{abbrvnat}

{
\appendix
\section{Appendix}
Next, we describe the remaining multi-period models proposed by \cite{Gun82}.

\begin{itemize}
	\item DSCLP---Dynamic Set Covering Location Problem.
	
	This problem is a multi-period SCP that restricts each demand point $j\in J$ 
	to be covered in time period $t_j\in \cal T$. 
	Hence, a facility covering $j$ should be opened 
	either at the beginning of that period 
	or in a previous one.
	In this problem,
	facilities cannot be closed during the planning horizon. 
	Furthermore, at most one facility/equipment
	can be 
	operating in each location.
	
	In order to adapt (GMSCLP) to the model (DSCLP) introduced by \cite{Gun82}, variables $y$, $z^\prime$, $w$ and $v$ 
	are
	removed from the model and 
	$y_{it}$ is replaced by $\sum_{\tau=1}^t z_{i\tau}$. 
	When performing these changes, constraints \eqref{yz}, \eqref{yz2}, \eqref{ww1}--\eqref{ybound}, \eqref{zpbound}--\eqref{vbound} are not longer needed.
	Besides, constraints \eqref{cover} are replaced by \eqref{ineqcov}. 
	Then, we set 
	$b_{jt} = 0$ for $j\in J$, $t < t_j$, where $t_j$ is the period in which demand of customer $j$ arises first. 
	Additionally, we take $b_{jt_j} = 1$ for $j\in J$ and, since facilities cannot be closed, we can
	set $b_{jt} = 0$ for $t > t_j$.
	The remaining parameters  from the model become:
	$c_{it}=f_{it}=0$, $o_{it}=o_t$ for $i\in I$, $t\in {\cal T}$ such that $o_{t}\geq o_{t+1}$ for $t=1,\ldots,T-1$.   
	In addition we set $p_t=m$ for $t\in {\cal T}$ and $e_i=1$ for $i\in I$.
	
	\item DSCLP2---Dynamic Set Covering Location Problem 2.
	
	This 
	is a multi-period set covering 
	problem where facilities cannot be closed during the planning horizon and at most 
	one facility can be opened  at each location. 
	In this problem, in each period $t\in {\cal T}$ a subset of demand points $J_t\subset J$ must be covered by at least one center. 
	Converting (GMSCLP) into (DSCLP2) is accomplished by setting
	$K_{jt}=K^\prime_{jt}=\emptyset$, $o_{it}=0$, $c_{it}=M$, $f_{it}=1$, $p_t=m$, $e_{i}=1$,  $\bar{y}_{i0}=0$ for $i\in I$, $j\in J$, $t\in {\cal T}$. 
	Finally, $b_{jt}=1$ for $j\in J_t$, $t\in {\cal T}$ and $b_{jt}=0$ for $j\in J \setminus J_t$, $t\in {\cal T}$. 
	
	Constraints \eqref{cover} are replaced by \eqref{ineqcov}. 
	Finally, since $K_{jt}$ and $K^\prime_{jt}$ are empty, constraints \eqref{ww1}--\eqref{ww3} are removed from the model.
	
	\item DSCPP---Dynamic Set Covering Phase Out Problem.
	
	In this problem, at the beginning of the planning horizon, there is one facility opened in each location $i\in I$. Using the terminology associated with (GMSCLP) we have $\bar{y}_{i0}=1$ for $i\in I$. 
	During the planning horizon, facilities can only be closed. 
	Similarly to (DSCLP2), a subset of customers $J_t\subset J$ must be covered in every period $t\in {\cal T}$. 
	In order to specialize (GMSCLP) to this model 
	we should consider only variables $y_{it}$ for $i\in I$, $t\in {\cal T}$ and we should set the parameters as follows: $K_{jt}=K^\prime_{jt}=\emptyset$, $o_{it}=M$, $c_{it}=0$, $f_{it}=1$, $p_t=m$, $e_i=1$, for $i\in I$, $j\in J$, $t\in {\cal T}$. 
	As in the previous model, $b_{jt}=1$ for $j\in J_t$, $t\in {\cal T}$ and $b_{jt}=0$ for $j\notin J_t$, $t\in {\cal T}$.
	
	\item DMCLP1---Dynamic Maximum Covering Location.
	
	This problem is a multi-period maximum covering location problem that seeks to minimize the weighted sum of demand points that are not covered in each time period. 
	In this case, we consider a set of empty locations ($I^o$) where a facility can be opened during the planning horizon and a set of locations ($I^c$) where a facility is currently open and that can be closed during the planning horizon.
	These two sets do not overlap.
	
	To model this problem using our general formulation we set $K_{jt}=\{1\}$.  Associated with each demand point, a population size $\varsigma_{jt}$ is considered and the benefit for covering this population is $g_{j1t}=-\varsigma_{jt}$.
	Besides, only $y$- and $w$-variables are considered.  
	Constraints \eqref{cover} are replaced by
	\begin{equation}
	\sum_{i\in I} a_{ijt}y_{it}\geq w_{j1t},\ \ \ j\in J,t\in {\cal T}. \label{ineqcov2}
	\end{equation}
	In this case, $f_{it}=0$, $h_{jt}=0$, $p_t=p_t$, $K'_{jt}=\emptyset$ and $e_i=1$ for $i\in I$, $j\in J$, $t\in {\cal T}$. 
	As before,  $o_{it}=M$ for $i \in I^c$ and $o_{it}$=0 otherwise for $t\in T$. Similarly, $c_{it}=M$ for $i\in I^o$ and $c_{it}=0$ otherwise for $t\in \cal{T}$.
	To obtain the same objective value as in the original model, the constant $-\sum_{j\in J}\sum_{t\in {\cal T}}g_{j1t}$ must be added to the objective value.
\end{itemize}

In tables \ref{const} and \ref{param}, we summarize the constraints and parameters, respectively, used to adapt (GMSCLP) to the models appearing in \cite{Gun82}.

\begin{table}
	\centering
	\caption{Constraints needed in analyzed formulations.}
	\begin{tabular}{lcl}
		& & \\[-2.0ex]
		\toprule
		Model  & & Constraints\\
		\midrule
		(COV)    & & 
		\eqref{sumy}, \eqref{cover},     \eqref{ybound},      \eqref{wbound}          \\
		(DSCLP)  & & \eqref{zbound}, \eqref{ineqcov}                                                  \\
		(DSCLP2), (DSCPP), (GDSCLP)  & & \eqref{yz},\eqref{yz2}, \eqref{ybound}-\eqref{zpbound}, \eqref{ineqcov}                               \\
		(DMCLP1), (DMCLP2) & & \eqref{sumy}--\eqref{yz2}, \eqref{ybound}--\eqref{wbound}, \eqref{ineqcov2} \\
		\bottomrule
	\end{tabular}%
	\label{const}
\end{table}%

\begin{sidewaystable}[htbp]
	\centering
	\caption{Paramaters to adapt model (GMSCLP) to existing covering location problems.}
	\scalebox{0.9}{
		\begin{tabular}{lllllllllllll}
			& & \\[-2.0ex]
			\toprule
			Model &$T$&$S$&$|K_{jt}^s|$&$|K_{jt}^{'s}|$&$c$&$f$&$g$&$p$&$b$&$e$&o\\
			\midrule
			(COV)&$1$&$1$&$K$&$0$&$0$&$f$&$g$&$p$&$b$&$e$&$0$\\
			(DSCLP)&$T$&$1$&$0$&$0$&$0$&$0$&$0$&$m$&$b_{jt_j}=1\textrm{ and }b_{jt}=0\textrm{ for }t\neq t_j,j\in J$&$1$&$o_{it}=o_t,o_t>o_{t+1},t\in {\cal T}$\\
			(DSCLP2)&$T$&$1$&$0$&$0$&$M$&$1$&$0$&$m$&$b_{jt}=1,t\in {\cal T},j\in J_t\textrm{ and }b_{jt}=0,t\in {\cal T},j\notin J_t$&$1$&$0$\\
			(DSCPP)&$T$&$1$&$0$&$0$&$0$&$1$&$0$&$m$&$b_{jt}=1,t\in {\cal T},j\in J_t\textrm{ and }b_{jt}=0,t\in {\cal T},j\notin J_t$&$1$&$M$\\
			(GDSCLP)&$T$&$1$&$0$&$0$&$c_{it}=M,\,i\in I^{o}$, $c_{it}=0\;i\in I^{c}$ &$1$&$0$&$m$&$b_{jt}=1,t\in {\cal T},j\in J_t\textrm{ and }b_{jt}=0,t\in {\cal T},j\notin J_t$&$1$&$o_{it}=M,\,i\in I^{c}$, $o_{it}=0\;i\in I^{o}$\\
			(DMCLP1)&$T$&$1$&$1$&$0$&$c_{it}=M,\,i\in I^{o}$, $c_{it}=0\;i\in I^{c}$&$0$&$-\varsigma$&$p$&$0$&$1$&$o_{it}=M,\,i\in I^{c}$, $o_{it}=0\;i\in I^{o}$\\
			(DMCLP2)&$T$&$1$&$1$&$0$&$c$&$0$& $-\varsigma$&$p$&$0$&$1$&$o$\\
			\bottomrule
		\end{tabular}%
	}
	\label{param}
\end{sidewaystable}%
}

\end{document}